\documentclass[pdflatex,sn-mathphys-num]{sn-jnl}

\usepackage[T1]{fontenc}
\usepackage{amsfonts,bm}
\usepackage{amssymb,amsmath}
\usepackage{amsthm}
\usepackage{graphicx,tikz}
\usepackage{mathrsfs}
\usepackage{array}
\usepackage{enumitem}
\usepackage{longtable}
\usepackage{url}
\usepackage{microtype}
\usepackage{needspace}

\renewcommand{\epsilon}{\varepsilon}
\renewcommand{\phi}{\varphi}

\renewcommand{\d}{\,\mathrm{d}}
\newcommand{\id}{\mathop{\bf id}}

\newtheoremstyle{thmstyleone}
{12pt plus2pt minus1pt}
{12pt plus2pt minus1pt}
{\itshape}
{0pt}
{\bfseries}
{.}
{.5em}
{}
\newtheoremstyle{thmstylethree}
{12pt plus2pt minus1pt}
{12pt plus2pt minus1pt}
{\normalfont}
{0pt}
{\bfseries}
{.}
{.5em}
{}

\theoremstyle{thmstyleone}
\ifcsname proposition\endcsname\else\newtheorem{proposition}{Proposition}\fi
\ifcsname theorem\endcsname\else\newtheorem{theorem}{Theorem}\fi
\ifcsname corollary\endcsname\else\fi
\ifcsname lemma\endcsname\else\newtheorem{lemma}{Lemma}\fi
\theoremstyle{thmstylethree}
\ifcsname assumption\endcsname\else\fi
\ifcsname definition\endcsname\else\newtheorem{definition}{Definition}\fi
\ifcsname remark\endcsname\else\newtheorem{remark}{Remark}\fi
\ifcsname example\endcsname\else\fi

\makeatletter
\AtBeginDocument{%
}
\makeatother

\newcommand{\spt}{\mathop{\rm spt}}

\renewcommand{\P}{\mathcal{P}}
\newcommand{\R}{\mathbb{R}}
\newcommand{\C}{\mathbb{C}}

\begin{document}

\title[Local Projected Controllability  of Nonlocal Continuity Equations]{\centering Local Projected Controllability\\ of Nonlocal Continuity Equations}

\author[1]{\fnm{Elena} \sur{Goncharova}}\email{goncha@icc.ru}
\author[1]{\fnm{Nikolay} \sur{Pogodaev}}\email{nickpogo@gmail.com}
\author*[2]{\fnm{Maxim} \sur{Staritsyn}}\email{staritsyn@fe.up.pt}

\affil[1]{\orgdiv{Matrosov Institute for System Dynamics and Control Theory}, \orgname{Siberian Branch of the Russian Academy of Sciences}, \orgaddress{\city{Irkutsk}, \country{Russia}}}
\affil[2]{\orgdiv{Research Center for Systems and Technologies (SYSTEC), ARISE, Department of Electrical and Computer Engineering}, \orgname{Faculdade de Engenharia, Universidade do Porto}, \orgaddress{\street{Rua Dr. Roberto Frias, s/n}, \postcode{4200-465}, \city{Porto}, \country{Portugal}}}

\abstract{We formulate and study the problem of local projected controllability (LPC) for characteristic solutions of controlled nonlocal continuity (transport) equations on the space of probability measures over a separable Banach space. The task is motivated in particular by self-interacting models in machine learning, whose action can be represented as controlled dynamics of the probability distribution over input data, while the actual output is given by a finite-dimensional vector of features extracted from the terminal distribution by a possibly nonlinear observation map. 
Given a fixed finite family of velocity fields which may depend on the transported measure, we ask whether the reachable values of the observable vector fill a neighborhood of its reference value. Our approach is to lift the transport equation in the space of measures to an ordinary differential equation on a Banach space of maps. This representation makes it possible to employ the machinery of Lie brackets developed for Banach manifolds, while the resulting sufficient condition for LPC can be formulated directly on the space of measures.
}

\keywords{nonlocal continuity equation, projected controllability, Lie brackets, dynamics of probability measures, lifting to a Banach space of maps}

\pacs[MSC Classification]{93B05, 93B27, 35F25, 28A33, 46T05}

\maketitle

\section{Introduction}

Classical controllability theory asks whether the system at hand has enough degrees of freedom to produce all desired motions \cite{agrachevControlTheoryGeometric2004} --- a property which is difficult to verify directly. This theory becomes much more complicated when the state space is infinite-dimensional. At the same time, in practice, 
one is often interested not in the entire infinite-dimensional state, but in a few numerical characteristics extracted from it. A question which is natural in many applications is the following: \emph{given a control system on a Banach manifold \(\mathcal M\), with control inputs $u$ and the corresponding states \(x^u_t\), along with an ``observation map'' \(g\colon \mathcal M \to \R^N\), \(N \in \mathbb N\), does the set \(\{g(x^u_t)\}\) of reachable observations contain a neighborhood of a fixed reference value?}

A problem of this type 
was addressed in \cite{Hermes1979}: the control system was given by an ODE in a Banach space, with a drift and a single control-linear term, and the observation was represented by a continuous {linear} map. The authors called this property \emph{local controllability of the observation} along the reference solution and gave a sufficient condition in terms of Lie directions projected by the map $g$. 
Subsequent results remained largely model-specific; various versions of the same task were studied, in particular, for the Euler and Navier--Stokes equations \cite{AgrachevSarychev2006,Shirikyan2007}.

In the present paper, a similar question is posed for a nonlocal transport equation on the space \(\P=\mathcal P(Q)\) of probability measures over a separable Banach space \(Q\). The state is now a measure \(\mu_t\in\P\) whose evolution is governed by a velocity field that may depend on the current measure and on the control; moreover, the observed quantity \(g\) may now be nonlinear. To simplify the terminology, we call the desired property \emph{local projected controllability} (LPC).\footnote{In that context, the word ``local'' has its usual topological meaning: the value of \(g\) is controlled in a neighborhood of the point \(g(\mu_0)\). This should not be confused with the distinction between local and nonlocal vector fields used in the theory of transport equations: a local vector field depends only on the point \(q\), whereas a nonlocal vector field may also depend on the current measure \(\mu\).}

\subsection{\texorpdfstring{{Background and related results}}{Background and related results}}

\bmhead{Classical controllability problem in finite and infinite dimensions}
A central result of classical geometric control theory is the Rashevskii--Chow theorem  \cite{Rashevskii1938,Chow1939}, which relates controllability of a driftless system to the properties of the Lie algebra generated by the controlled vector fields. In $\R^n$, one can choose finitely many Lie brackets whose values form a basis of the tangent space and then apply the inverse function theorem to a finite composition of the corresponding flows; see, e.g.,~\cite{agrachevControlTheoryGeometric2004}. 

A passage to infinite dimensions changes the situation substantially: the resulting Lie algebra may be dense in the tangent space without spanning it, and this is no longer sufficient for the reachable set to contain a neighborhood of the initial state. One of the first infinite-dimensional counterparts of the classical Rashevskii--Chow theorem was obtained by~\cite{DudnikovSamborski1980} for systems in Banach spaces; related sufficient conditions for global approximate controllability were established in~\cite{Ledyaev2004}. Further extensions were made for control systems on Banach manifolds~\cite{KipkaLedyaev2015} and regular manifolds modeled on convenient locally convex spaces~\cite{KhajehSalehaniMarkina2014}. 
Related Lie-algebraic controllability results for ensembles of points and families of mappings can be found in~\cite{agrachevControlSpacesEnsembles2020,agrachevControlManifoldsMappings2022}.

\bmhead{Dynamics in spaces of measures}
The evolution of a probability measure under the action of a vector field admits several complementary descriptions emphasizing different aspects of the same dynamics.  In the Eulerian formulation, the measure itself is the unknown and its evolution is governed by a continuity equation. The Lagrangian formulation focuses on characteristic curves and reconstructs the measure as the push-forward of its initial value. The same evolution can also be considered intrinsically as a curve in a space of probability measures endowed with a Kantorovich metric; the relations between these formulations are studied in~\cite{cavagnariLagrangianEulerianKantorovich2022,AVERBOUKH2025113205}. In this paper, we represent the dynamics on $\mathcal P$ as an ODE on the Banach space $C_b(Q;Q)$ of bounded continuous maps from $Q$ to itself.

\bmhead{Differential structures on spaces of measures}
Our lift of the dynamics on $\mathcal P$ is motivated, in particular, by the fact that $\mathcal P$ and its subspaces lack natural structures of a Banach manifold, which prevents a direct application of the geometric control results described above. A promising candidate for establishing such a structure would be the subspace $\mathcal P_2(\R^n)$ of probability measures on $\R^n$ with finite second moment, endowed with the quadratic Kantorovich (Wasserstein) distance $W_2$: a formal Riemannian calculus on this space was developed in~\cite{ottoGeometryDissipativeEvolution2001} (see also \cite{AGS}). However, in this calculus, tangent velocities at a point $\mu$ are represented by elements of a Hilbert space contained in $L^2(\mu;\R^n)$, which changes with the base measure; furthermore, the metric tangent cone may be nonlinear at singular measures; see, e.g., \cite{lottTangentConesWasserstein2016}.

\bmhead{Controllability for nonlocal continuity equations}
The controllability of continuity equations and, more generally, of evolution equations on spaces of measures has been considered from several points of view. Approximate and exact controllability for a continuity equation with a localized control vector field are proved in~\cite{duprezApproximateExactControllability2019}, and the corresponding minimum-time problem is treated in~\cite{duprezMinimalTime2020}. For nonlocal continuity equations, stabilization of trajectories by localized controls is considered in~\cite{PogodaevRossi2024TrajectoryStabilization}, while stabilization by geodesic perturbations is obtained in~\cite{PogodaevRossi2026GeodesicStabilization}. Attainability of prescribed sets of probability measures is studied in~\cite{CavagnariMarigonda2021Attainability}. 
Continuity equations also arise naturally in control problems for ensembles of initial states; see, for instance,~\cite{StaritsynPogodaevChertovskihPereira2021,ChertovskihPogodaevStaritsynAguiar2023}.
Related control problems in Kantorovich (Wasserstein) spaces are treated in~\cite{BonnetFrankowska2021a,JimenezMarigondaQuincampoix2020}.
To the best of our knowledge, problems of the LPC-type have never been addressed in this context.

\bmhead{Mean-field models of Transformers}
A practical motivation for this study comes from self-attention in Transformers~--- a machine-learning architecture introduced in~\cite{Vaswani2017}: a finite collection of vectors in $\R^n$ is naturally represented by its empirical measure, subject to nonlocal transformations on $\mathcal P(\R^n)$ via so-called self-attention maps. A terminal output map extracts a finite-dimensional feature vector from the resulting measure, and controlling this vector leads precisely to the LPC problem studied in this paper (see \S~\ref{sec:transformers}). Related problems of controllability, interpolation, and clustering are considered in~\cite{10540044,Agrachev2025,GeshkovskiRigolletRuizBalet2026,AlcaldeFantuzziZuazua2026,Cuchiero}. 

\bmhead{Symmetry restrictions}
Symmetry constraints are customarily imposed on the dynamics in order to reflect physical, constructive or any other internal properties of the modeled object; systems respecting such constraints are termed \emph{equivariant}. Controllability of equivariant systems is studied in~\cite{Agrachev2025} with applications to particle systems and artificial neural networks; in the latter context see also~\cite{Puny2022,Kaba2023}. In Section~\ref{sec:equivariance}, we address a related but different question: local and nonlocal vector fields, in the sense of the theory of transport equations, are considered under the same symmetry constraints. 

\subsection{\texorpdfstring{{Main results}}{Main results}}

Our principal result is a sufficient condition for LPC of characteristic solutions to a controlled nonlocal continuity equation. We obtain it by lifting the transport equation to an ODE on the Banach space of displacement maps. From the solution of the lifted equation, we reconstruct the characteristics and the corresponding curve of measures. Theorem~\ref{theorem} extends the sufficient LPC criterion of~\cite{Hermes1979} to driftless ODEs with vector-valued controls on Banach manifolds.

For a nonlocal vector field (NVF) $V$, we consider its action $L_V$ on functionals of a measure. The principal criterion, Theorem~\ref{thm:transport-projected-rank}, is then stated directly on $\mathcal P$: the values obtained by applying iterated commutators of the operators $L_{V_j}$ to the observable at the reference measure must span $\R^N$. Smoothness is imposed locally on the lifts of the NVFs and observables; for integral functionals and NVFs with bounded kernels, it follows from Lemma~\ref{lem:smooth-class}. If a finite family of functionals is closed under the operators $L_{V_j}$, the corresponding projected dynamics admits an exact finite-dimensional reduction. For integral functionals and NVFs defined by interaction kernels, the operators $L_V$ and their commutators admit explicit formulas in terms of those kernels.

The main results are illustrated by three examples. In the first one, the missing projected direction appears at the first commutator in a surrogate model of a Transformer. In the other two, under translation and rotation equivariance, respectively, the observable is fixed by every admissible local system, whereas suitable NVFs attain the required projected rank.

\subsection{\texorpdfstring{{Organization of the paper}}{Organization of the paper}}

The rest of the paper is organized as follows. In  Section~\ref{sec:prelim}, we formulate the LPC problem on a Banach manifold and provide a general sufficient condition. Section~\ref{sec:rep} is devoted to the dynamics on the space of measures given by nonlocal transport equations. Section~\ref{sec:main-results} contains a detailed exposition of the announced lifting of the transport equation to an ODE on a Banach space, followed by the corresponding calculus of Lie brackets. Section~\ref{sec:transformers} relates this framework to context-dependent models of Transformers. This line is completed by a discussion of the LPC under extra symmetry constraints in Section~\ref{sec:equivariance}. The proofs of the principal technical statements are collected in the Appendix.

\section{Preliminaries}\label{sec:prelim}

Given Banach spaces $Y,Z$, we write $\mathcal L(Y;Z)$ for the space of bounded linear maps with the operator norm. A map $J$ from an open set $\mathcal V\subset Y$ to $Z$ is termed \emph{smooth} if it has continuous Fr\'echet derivatives of all orders. The derivative of $J$ at $e\in\mathcal V$ is $DJ(e)\in\mathcal L(Y;Z)$, whose value at $v\in Y$ is denoted by $DJ(e)[v]$.

The terminology of smooth Banach manifolds, tangent vectors, vector fields (VFs) and Lie brackets follows~\cite[pp.~22--154]{Lang1999}; the relevant conventions are recalled in Appendix~\ref{app:manifold-notation}. The tangent bundle of $\mathcal M$ is denoted by $\mathrm T\mathcal M$. For a smooth scalar function $\phi$ and a VF $F$, the notation $F\phi$ stands for the derivative of $\phi$ along $F$. The Lie bracket is taken with the convention $[F_1,F_2]_{\mathcal M}\phi=F_1(F_2\phi)-F_2(F_1\phi)$.

\subsection{Local projected controllability (LPC)}

Let $m\geq1$ be an integer. Given a real smooth Banach manifold
$\mathcal M$, a compact set $U\subset\R^m$ containing the cube
$[-1,1]^m$, and smooth VFs $v_j\colon\mathcal M\to\mathrm T\mathcal M$,
$j=1,\ldots,m$, we define the admissible controls by $\mathcal U_T\doteq L^\infty([0,T];U)$, $T>0$, and consider the control-linear system
\begin{equation}
\dot x(t)
=
\sum_{j=1}^m u_j(t)\, v_j\bigl(x(t)\bigr),
\qquad
u\in\mathcal U_T.
\label{ODE}
\end{equation}

Fix $x_0\in\mathcal M$. Solutions of~\eqref{ODE} are understood in the Carath\'eodory sense: in a chart $(\mathcal O,\psi)$, their coordinate representations $\xi=\psi\circ x$ are absolutely continuous and satisfy the corresponding Bochner integral equation, whose integrand is strongly measurable since the controls are measurable and bounded and the coordinate representations of the VFs are continuous.

Since the coordinate representations of the VFs are locally bounded and locally Lipschitz, and $U$ is compact, the existence and uniqueness theorem for ODEs in Banach spaces~\cite[pp.~4--17]{Deimling1977} gives a chart $(\mathcal O,\psi)$ containing $x_0$ and a number $T_0>0$ such that, for every $T\in(0,T_0]$ and $u\in\mathcal U_T$, system~\eqref{ODE} has a unique solution $x^u(\cdot;x_0)\colon[0,T]\to\mathcal O$, with $x^u(0;x_0)=x_0$. In particular, $T_0$ can be chosen independently of $u$.

Let $N\geq1$ be an integer, and
$g\colon\mathcal M\to\R^N$ a smooth map. 
\begin{definition}
We say that the ODE \eqref{ODE} is \emph{locally projected
controllable at a point $x_0\in\mathcal M$} (w.r.t. $g$) if there exist \(T>0\) and a neighborhood
\(\mathcal N\) of the point \(g(x_0)\) in \(\R^N\) such that, for every
\(y_T\in\mathcal N\), there exists a control \(u\in\mathcal U_T\) for which the corresponding solution is defined on $[0,T]$ and $g\bigl(x^u(T;x_0)\bigr)=y_T$.
\end{definition}

Denote by $\mathrm{Lie}\{v_1,\dots,v_m\}$ the Lie algebra generated by
$v_1,\ldots,v_m$.

\begin{theorem}[A sufficient condition for LPC]\label{theorem}
Suppose the Lie algebra $\mathrm{Lie}\{v_1,\ldots,v_m\}$ contains VFs $W_1,\ldots,W_N$ satisfying the following projected rank condition:
\begin{equation}\label{c:pr}
\operatorname{span}\bigl\{
Dg(x_0)[W_i(x_0)]\colon i=1,\ldots,N
\bigr\}
=\R^N.
\end{equation}
Then, for every sufficiently small $T>0$, the values $g\bigl(x^u(T;x_0)\bigr)$ attained by piecewise-constant controls contain a neighborhood of $g(x_0)$. In particular, system~\eqref{ODE} is LPC at $x_0$ w.r.t.~$g$.
\end{theorem}

\begin{remark}[Computational price of commutators]
No universal upper bound on the length of the brackets needed in~\eqref{c:pr} is known, and testing long brackets may be difficult. A commutator loop of two VFs has duration of order $\varepsilon$ but produces a bracket displacement of order $\varepsilon^2$; for a nested bracket containing $k$ VFs, the first nonvanishing contribution is of order $\varepsilon^k$. Higher brackets thus require more switchings and smaller displacements on the same time scale. Switching restrictions, costs and numerical errors may limit the realizable depth.
\end{remark}

\begin{remark}[Quantitative characterization of LPC]\label{re:numerical}
We have already stressed that controllability is a qualitative property. One possible way to quantify LPC is to fix a family of scalar observables separating measures and count how many can be locally controlled simultaneously. The resulting number may depend on the chosen family, and controllability of arbitrarily large finite families does not imply an actual infinite-dimensional controllability, since the corresponding control horizons $T$ and neighborhoods need not be uniform. 
\end{remark}

\section{Dynamics on the space of measures}\label{sec:rep}

\subsection{Conventions}

Throughout, $Q$ is a separable real Banach space, and $\P\doteq\mathcal P(Q)$ is the set of Borel probability measures on $Q$. Let $E$ be a Banach space. For $k\geq1$ and a Borel map $\phi\colon Q^k\to E$, put
\[
\langle\phi\rangle(\mu)\doteq\int_{Q^k}\phi\d\mu^{\otimes k},
\]
whenever the integral exists. The number $k$ is understood from the domain of the kernel $\phi$ and is therefore omitted from the notation. This notation will be convenient when we compute derivatives of polynomial functionals in Section~\ref{ssec:computing-rank}.

We denote the identity map on $Q$ by ${\id}_Q$. Given a Borel map $S\colon Q\to Q$, the push-forward of $\mu$ by $S$ is the measure $S_\sharp\mu \in \mathcal P$ defined by the relation $\langle\phi\rangle(S_\sharp\mu)=\langle\phi\circ S\rangle(\mu)$ for all bounded Borel functions $\phi\colon Q\to\R$.

\Needspace{7\baselineskip}
\subsection{Nonlocal continuity equations and LPC}

By an NVF on $Q$ we mean a map $V\colon Q\times\P\to Q$; it is a local VF if it does not depend on the measure. Given NVFs $V_1,\ldots,V_m$ and a compact set $U\subset\R^m$ containing $[-1,1]^m$, we define the controlled NVF $F\colon Q\times\P\times U\to Q$ by
\begin{equation}\label{eq:controlled-nonlocal-field}
F(q,\mu,\mathrm u)\doteq\sum_{j=1}^m\mathrm u_jV_j(q,\mu),\qquad\mathrm u\in U.
\end{equation}
The NVF $F$ determines the formal partial differential equation
\begin{equation}\label{eq:CE}
\partial_t\mu_t+\mathrm{div}\bigl(F(\cdot,\mu_t,u(t))\,\mu_t\bigr)=0,
\qquad u\in\mathcal U_T,
\end{equation}
known as the \emph{nonlocal continuity equation} on $\mathcal P$.

Fix $\vartheta\in\P$ and $u\in\mathcal U_T$. A curve $t\mapsto\mu_t\in\P$ is called a \emph{distributional solution} of~\eqref{eq:CE} with initial value $\vartheta$ if $\mu_0=\vartheta$ and, for every smooth $\phi\colon Q\to\R$ bounded together with its derivative $D\phi$, the function $t\mapsto\langle\phi\rangle(\mu_t)$ is absolutely continuous and, for almost every (a.e.) $t\in[0,T]$,
\begin{equation}\label{eq:weak-formulation}
\frac{\d}{\d t}\int_Q\phi(q)\d\mu_t(q)
=\int_Q D\phi(q)[F(q,\mu_t,u(t))]\d\mu_t(q).
\end{equation}

An important role in what follows will be played by the Banach space \[\mathscr Q\doteq C_b(Q;Q)
\] of bounded continuous maps $Q\to Q$, endowed with the supremum norm $\|\cdot\|_\infty$. 

Let $t\mapsto\Xi_t^{u,\vartheta}(q)$ be a solution to the Bochner integral equation
\begin{equation}\label{eq:nlf}
\Xi_t^{u,\vartheta}(q)
=
q+
\int_0^t
F\bigl(\Xi_s^{u,\vartheta}(q),(\Xi_s^{u,\vartheta})_\sharp\vartheta,u(s)\bigr)\d s,
\qquad q\in Q.
\end{equation}
We call such maps the \emph{characteristics} of~\eqref{eq:CE} (issued from $q$); their local existence and uniqueness are established in Proposition~\ref{prop:local-lift}. 
The corresponding curve of measures
\begin{equation}\label{eflowsol}
t\mapsto\mu_t^{u,\vartheta}\in\P, \qquad  \mu_t^{u,\vartheta}\doteq(\Xi_t^{u,\vartheta})_\sharp\vartheta,
\end{equation}
is called a \emph{characteristic solution} of~\eqref{eq:CE} with initial value $\mu_0 = \vartheta$\footnote{Under standard Cauchy--Lipschitz conditions, the uniqueness of a distributional solution to~\eqref{eq:CE} on $\mathcal P(\R^n)$ follows from the superposition principle~\cite[Theorem~8.2.1]{AGS} and the uniqueness of characteristics. An extension of this argument to infinite-dimensional Banach spaces requires additional assumptions; see, e.g.,~\cite{STEPANOV20171044}.}; when $\vartheta$ is fixed, we write $\Xi_t^u$ and $\mu_t^u$, omitting the reference measure from the notation.

We now pose the LPC problem for~\eqref{eq:CE} in the class of characteristic solutions.
\begin{definition}
Let $g$ be a map with values in $\R^N$, defined on a set $\mathcal G\subset\P$ containing $\vartheta$. The equation~\eqref{eq:CE} is \emph{locally projected controllable at $\vartheta$ w.r.t.\ $g$ in the class of characteristic solutions} if there exist $T>0$ and a neighborhood $\mathcal N$ of $g(\vartheta)$ in $\R^N$ such that, for every $y\in\mathcal N$, one can find a control $u\in\mathcal U_T$ for which the characteristic solution~\eqref{eflowsol} is defined on $[0,T]$, $\mu_T^u\in\mathcal G$, and $g(\mu_T^u)=y$.
\end{definition}

\section{Main results}\label{sec:main-results}

\subsection{\texorpdfstring{{Lifting to an ODE on a Banach space}}{Lifting to an ODE on a Banach space}}

Fix the reference measure $\vartheta\in\P$. For $X\in\mathscr Q$, put
\[
\widetilde X\doteq{\id}_Q+X,
\qquad
\vartheta^X\doteq\widetilde X_\sharp\vartheta,
\]
and denote by $0_{\mathscr Q}$ the zero element of $\mathscr Q$. Consider a map $f\colon\mathcal G\to E$, with $\vartheta\in\mathcal G\subset\P$, and an NVF $V\colon Q\times\mathcal G\to Q$. Their \emph{lifts} at $\vartheta$ are introduced as 
\[
\widehat f\colon\mathcal D\to E, \quad \widehat f(X)\doteq f(\vartheta^X)\quad\mbox{ and }\ \quad
\widehat V\colon\mathcal D\to\mathscr Q, \quad \widehat V(X)\doteq\bigl(q\mapsto V(\widetilde X(q),\vartheta^X)\bigr),
\]
where $\mathcal D$ is any open neighborhood of $0_{\mathscr Q}$ on which $\vartheta^X\in\mathcal G$ and the second expression belongs to $\mathscr Q$. The same construction applies to $f$ and $V$ on the respective domains.

\begin{definition}[Smoothness]\label{def:smooth}
A map $f\colon\mathcal G\to E$ or an NVF $V\colon Q\times\mathcal G\to Q$ is termed \emph{smooth at $\vartheta$} if its lift is smooth on some open neighborhood of $0_{\mathscr Q}$. We say that $f$ or $V$ is \emph{smooth on $\mathcal G$} if it is smooth at all $\mu\in\mathcal G$.
\end{definition}

\begin{remark}
Definition~\ref{def:smooth} does not refer to any topology on $\P$: the smoothness of $f$ at $\vartheta$ is determined by its lift to $\mathscr Q$, rather than by its behavior on a neighborhood of $\vartheta$ in the narrow topology. It is nevertheless preserved under sufficiently small transport perturbations of $\vartheta$. Indeed, for sufficiently small $Y,Z\in\mathscr Q$, we have $({\id}_Q+Z)_\sharp\vartheta^Y=\vartheta^{Y+Z\circ\widetilde Y}$, and the map $Z\mapsto Y+Z\circ\widetilde Y$ is affine and continuous on $\mathscr Q$. Thus $f$ is smooth at $\vartheta^Y$ for all sufficiently small $Y\in\mathscr Q$. The same conclusion does not follow for an NVF from the smoothness of its lift alone: changing the reference measure may require values of the NVF at points outside the range of $\widetilde Y$.
\end{remark}

\Needspace{10\baselineskip}
\begin{lemma}\label{lem:smooth-class}
Let $\phi\colon Q^{1+k}\to E$ be smooth and bounded together with its derivatives of all orders. Then the map from $\mathscr Q$ to $C_b(Q;E)$, defined by
\[
X\longmapsto\Bigl[q\mapsto\int_{Q^k}\phi\bigl(\widetilde X(q),\widetilde X(r_1),\ldots,\widetilde X(r_k)\bigr)\d\vartheta^{\otimes k}(r_1,\ldots,r_k)\Bigr],
\]
is smooth. In particular, integral maps $\mu\mapsto\langle\phi\rangle(\mu)$ with kernels of this class, as well as NVFs defined by such kernels, are smooth. Smoothness is preserved by sums and, for scalar-valued functionals, by products and quotients wherever the denominator stays away from zero. The same holds for NVFs of the form $(q,\mu)\mapsto W(q,f_1(\mu),\ldots,f_K(\mu))$, where $f_i\colon\mathcal G\to\R$ are smooth, while $W\colon Q\times\R^K\to Q$ is smooth and bounded together with its derivatives.
\end{lemma}

On a common open neighborhood $\mathcal D$ of $0_{\mathscr Q}$, we consider the control system generated by the lifted VFs $\widehat V_1,\ldots,\widehat V_m$:
\begin{equation}\label{eq:lifted-control-system}
\dot X_t
=
\sum_{j=1}^m u_j(t)\widehat V_j(X_t),
\qquad X_0=0_{\mathscr Q}.
\end{equation}
We interpret~\eqref{eq:lifted-control-system} as a Bochner integral equation in $\mathscr Q$.
\begin{proposition}\label{prop:local-lift}
Suppose that $V_1,\ldots,V_m$ are smooth at $\vartheta$. Then there exists $T_0>0$ such that, for all $T\in(0,T_0]$ and $u\in\mathcal U_T$, system~\eqref{eq:lifted-control-system} has a unique solution $X^u\in C([0,T];\mathcal D)$. Define
\begin{equation}\label{eq:lift-flow-reconstruction}
\Xi_t^{u,\vartheta}\doteq{\id}_Q+X_t^u,
\qquad
\mu_t^{u,\vartheta}\doteq\vartheta^{X_t^u}.
\end{equation}
The map $t \mapsto \Xi^{u,\vartheta}_t$ is the unique solution of the characteristic equation~\eqref{eq:nlf} satisfying the inclusion $(\Xi^{u,\vartheta}-{\id}_Q)\in C([0,T];\mathscr Q)$. The curve $t\mapsto\mu_t^{u,\vartheta}$ is a distributional solution of~\eqref{eq:CE}.
\end{proposition}

\begin{remark}[Relation to Lions' lift]
The idea of lifting a nonlinear object to a suitable linear space is quite standard. Our realization of this idea mimics, in essence, that of P.-L.~Lions, who developed a calculus on $\mathcal P_2(\R^n)$ by lifting a functional of a probability measure to the Hilbert space of square-integrable random variables through the law map; see, e.g.,~\cite[Chapter~5]{CarmonaDelarue2018}. In our approach, the lift is applied to the characteristic dynamics itself: the reference measure is fixed, while the bounded displacement $X$ parametrizes both the characteristic map and the transported measure.
\end{remark}

\subsection{Derivatives along nonlocal vector fields and the rank condition}

Given a map $f\colon\mathcal G\to E$ and an NVF $V$, we define, whenever the expression exists, the derivative of $f$ along $V$ by the equality
\begin{equation*}
(L_Vf)(\mu)
\doteq
\left.\frac{\d}{\d\varepsilon}\right|_{\varepsilon=0}
f\!\left(
[{\id}_Q+\varepsilon V(\cdot,\mu)]_\sharp\mu
\right).
\end{equation*}
For $f$ and $V$ smooth at $\vartheta$, the identity $[{\id}_Q+\varepsilon V(\cdot,\vartheta^X)]\circ\widetilde X=\widetilde X+\varepsilon\widehat V(X)$ and the chain rule yield
\begin{equation}\label{eq:derivation-lift}
\widehat{L_Vf}(X)
=
D\widehat{f}(X)[\widehat V(X)],
\qquad X\in\mathcal D.
\end{equation}
Thus the lift of $L_Vf$ is the directional derivative of $\widehat f$ along the VF $\widehat V$ on the Banach space $\mathscr Q$. In particular, $L_Vf$ is again smooth at $\vartheta$.

\begin{lemma}\label{lem:LV-chain-rule}
Suppose that $V_1,\ldots,V_m$ and $f$ are smooth at $\vartheta$. Then, for all sufficiently small $T>0$ and any $u\in\mathcal U_T$, the function $t\mapsto f(\mu_t^u)$ is absolutely continuous and, for a.e.~$t$,
\begin{equation*}
\frac{\d}{\d t}f(\mu_t^u)
=
\sum_{j=1}^m u_j(t)(L_jf)(\mu_t^u),
\qquad L_j\doteq L_{V_j}.
\end{equation*}
\end{lemma}

On the class of functionals smooth at $\vartheta$, we consider commutators $[A,B]\doteq AB-BA$ of the operators $L_j$. Denote by $\operatorname{Lie}\{L_1,\ldots,L_m\}$ the Lie algebra they generate. By~\eqref{eq:derivation-lift}, each $L_j$ lifts to the directional derivative along $\widehat V_j$. Since the Lie bracket of VFs on a Banach space is the commutator of the corresponding directional derivatives, every combination $P$ of iterated commutators of the $L_j$ lifts to the directional derivative along the same combination $\widehat P$ of iterated Lie brackets of the $\widehat V_j$:
\begin{equation}\label{eq:Lie-derivation-lift}
\widehat{Pf}(X)
=
D\widehat{f}(X)[\widehat P(X)].
\end{equation}

For NVFs $V,W$ on $\R^n$, let
\begin{equation}\label{eq:material}
(\nabla_VW)(q,\mu)\doteq\left.\frac{\d}{\d\varepsilon}\right|_{\varepsilon=0}
W\Bigl(q+\varepsilon V(q,\mu),\bigl[{\id}_Q+\varepsilon V(\cdot,\mu)\bigr]_\sharp\mu\Bigr)
\end{equation}
be the rate of change of $W$ when the point $q$ and the whole measure $\mu$ move along $V$.

\begin{proposition}[Lie bracket of NVFs]\label{prop:bracket}
Let $Q=\R^n$, and let $V,W$ be smooth at $\vartheta$. Then $[\widehat V,\widehat W]$ is the lift of the NVF $[V,W]\doteq\nabla_VW-\nabla_WV$.
\end{proposition}

Together with~\eqref{eq:derivation-lift}, this gives $[L_V,L_W]f=L_{[V,W]}f$ at the measures $\vartheta^X$, $X\in\mathcal D$.

\begin{remark}[Intrinsic derivatives]
Let $Q=\R^n$. Suppose that an NVF $W$ admits a variational derivative w.r.t.\ the measure in the sense of~\cite{kolokoltsov2019,CardMaster2019}, i.e., a function $\frac{\delta W}{\delta\mu}(q,\mu)(r)$ such that
\[
\left.\frac{\d}{\d\varepsilon}\right|_{\varepsilon=0+}
W(q,(1-\varepsilon)\mu+\varepsilon\nu)
=
\int_{\R^n}\frac{\delta W}{\delta\mu}(q,\mu)(r)\d(\nu-\mu)(r),\qquad\nu\in\P,
\]
and let $\bm D_\mu W(q,\mu)(r)\doteq D_r\frac{\delta W}{\delta\mu}(q,\mu)(r)$ be its intrinsic derivative. If $D_qW$ and $\bm D_\mu W$ exist and are regular enough for the chain rule along transport perturbations of $\mu$ (e.g., for NVFs defined by smooth bounded kernels), then
\[
\nabla_VW(q,\mu)
=
D_qW(q,\mu)[V(q,\mu)]
+\int_{\R^n}\bm D_\mu W(q,\mu)(r)[V(r,\mu)]\d\mu(r),
\]
and therefore
\[
\begin{aligned}
[V,W]
&=
D_qW[V]-D_qV[W]\\
&+
\int_{\R^n}\bigl(
\bm D_\mu W(q,\mu)(r)[V(r,\mu)]
-\bm D_\mu V(q,\mu)(r)[W(r,\mu)]
\bigr)\d\mu(r),
\end{aligned}
\]
where $V,W,D_qV,D_qW$ are evaluated at $(q,\mu)$. The first two terms form the usual Lie bracket of VFs on $\R^n$; the integral term accounts for the dependence on the measure and vanishes for local VFs.
\end{remark}

\begin{theorem}[A sufficient condition for LPC of nonlocal continuity equations]\label{thm:transport-projected-rank}
Let $V_1,\ldots,V_m$ and the components of $g$ be smooth at $\vartheta$. If there exist $P_1,\ldots,P_N\in\operatorname{Lie}\{L_1,\ldots,L_m\}$ such that
\begin{equation}\label{eq:transport-projected-rank}
\operatorname{span}\left\{(P_i g)(\vartheta)\colon i=1,\ldots,N\right\}=\R^N,
\end{equation}
then~\eqref{eq:CE} is locally projected controllable at $\vartheta$ w.r.t.\ $g$ in the class of characteristic solutions. The conclusion holds with piecewise-constant controls and arbitrarily small control time.
\end{theorem}

\subsection{Computing the projected directions}\label{ssec:computing-rank}

In general, $L_jg$ depends on the whole measure, not merely on $g$: an equation for the mean may involve the second moment, while an equation for the second moment may involve higher moments. This hierarchy may, however, close after finitely many steps. For example, on $\mathcal P_2(\R)$, put
\[
\mathbb E(\mu)\doteq\int_\R q\d\mu(q),
\qquad
\mathbb V(\mu)\doteq\int_\R(q-\mathbb E(\mu))^2\d\mu(q).
\]
Fix $\vartheta\in\mathcal P_2(\R)$. The NVF $V\colon\R\times\mathcal P_2(\R)\to\R$, defined by $V(q,\mu)\doteq\mathbb V(\mu)$, generates translations at a speed given by the variance. Since $L_V\mathbb E=\mathbb V$ and $L_V\mathbb V=0$, the pair $(\mathbb E,\mathbb V)$ satisfies the closed system $\dot{\mathbb E}=u\,\mathbb V$, $\dot{\mathbb V}=0$, whereas the mean alone does not evolve in closed form. Thus the mean is locally controllable if and only if $\mathbb V(\vartheta)>0$; Theorem~\ref{thm:transport-projected-rank} gives only the ``if'' part. More generally, starting from $g$, one may successively add the functionals appearing in the derivatives $L_jg$. If this procedure yields a finite family closed under all $L_j$, the following proposition applies.

\begin{proposition}[Closure of a finite family of functionals]\label{prop:closure}
Let $V_1,\ldots,V_m$ and $f=(f_1,\ldots,f_K)$ be smooth at $\vartheta$, let $h\colon\R^K\to\R^N$ be smooth, and put $g\doteq h\circ f$. Suppose that, on a common open neighborhood $\mathcal D$ of $0_{\mathscr Q}$ where the lifts are smooth, the following closure relations hold:
\begin{equation}\label{eq:closure-relation}
L_j f(\mu)=a_j(f(\mu)),
\qquad j=1,\ldots,m,
\end{equation}
for every $\mu=\vartheta^X$, $X\in\mathcal D$, where the maps $a_j\colon\R^K\to\R^K$ are smooth. Then, for every sufficiently small $T>0$ and every $u\in\mathcal U_T$, the function $y(t)\doteq f(\mu_t^u)$ is the unique solution of
\begin{equation}\label{eq:closed-system}
\dot y(t)=\sum_{j=1}^m u_j(t)a_j(y(t)),
\qquad y(0)=f(\vartheta).
\end{equation}
In particular, the values of $g(\mu_T^u)$ and of $h(y(T))$ reachable at each such time $T$ coincide.
\end{proposition}

\begin{proposition}[Integral calculus]\label{prop:polynomial-calculus}
Suppose that $\phi\colon Q^k\to E$ is smooth and bounded together with its derivatives. Given an NVF $V$ such that $V(\cdot,\mu)$ is bounded and Borel, we have
\begin{equation}\label{eq:polynomial-calculus}
(L_V\langle\phi\rangle)(\mu)
=
\sum_{a=1}^k
\int_{Q^k}
D_{q_a}\phi(q_1,\ldots,q_k)
\bigl[V(q_a,\mu)\bigr]
\d\mu^{\otimes k}.
\end{equation}
If, in addition, $V(q,\mu)=\langle K(q,\cdot)\rangle(\mu)$, where $K\colon Q\times Q\to Q$ is smooth and bounded together with its derivatives, it holds $L_V\langle\phi\rangle=\langle\partial_K\phi\rangle$, in which
\[
(\partial_K\phi)(q_1,\ldots,q_{k+1})
\doteq
\sum_{a=1}^k
D_{q_a}\phi(q_1,\ldots,q_k)
\bigl[K(q_a,q_{k+1})\bigr].
\]
In particular, if $V_i$ and $V_j$ correspond to kernels $K_i$ and $K_j$, then
\begin{equation}\label{eq:interaction-commutator}
[L_i,L_j]\langle\phi\rangle
=
\langle\partial_{K_i}\partial_{K_j}\phi-\partial_{K_j}\partial_{K_i}\phi\rangle,
\end{equation}
and each iterated commutator containing $p$ operators $L_j$ applied to $\langle\phi\rangle$ is an integral over $Q^{k+p}$. By the chain and Leibniz rules, the same formulas apply to smooth functions of finitely many integral functionals $\langle\phi_1\rangle,\ldots,\langle\phi_K\rangle$.
\end{proposition}

\section{A control-theoretic view of Transformers}
\label{sec:transformers}

A concrete case in which LPC on the space of measures comes into play is the next-token prediction problem in language models.
In such models, the text is viewed as a sequence of elementary units called \emph{tokens}, represented by vectors $q_1,\ldots,q_s$, $s\geq 1$, in a vector space $Q\doteq \R^n$ with $n$ typically very large. A paradigmatic task is to define the probability distribution of the next element $q_{s+1}$ based on the ``context'' --- the relationships between the elements of the input sequence $\{q_1,\ldots,q_s\}$.

A standard approach to this type of task is based on artificial (typically, residual) neural networks termed \emph{Transformers}. The operation of \emph{self-attention} evaluates the relative importance of the other tokens w.r.t. a given one~\cite{Vaswani2017}~--- via linear maps $\mathsf Q, \mathsf K\colon\R^n\to\R^d$ and $\mathsf V\colon\R^n\to\R^n$, called, in a mnemonic way, the \emph{query}, \emph{key} and \emph{value} maps.
Roughly speaking, the query $\mathsf Q q_i$ describes what a token represented by $q_i$ is ``looking for'' in the other ones, the key $\mathsf K q_j$ says what another token represented by $q_j$ has to offer in this comparison, and the value $\mathsf V q_j$ represents the information contributed by this token to the update of the $i$-th token.
The update of the $i$-th token takes the form
\(
q_i\longmapsto q_i+\mathfrak A_i(q_1,\ldots,q_s),\)
$\mathfrak A_i(q_1,\ldots,q_s)
=\sum_{j=1}^s\alpha_{ij}\mathsf V q_j,$
where the attention weights are defined by
\[
\alpha_{ij}=\frac{\exp\bigl(\langle\mathsf Qq_i,\mathsf Kq_j\rangle/\sqrt d\bigr)}{\mathfrak B_i},
\qquad
\mathfrak B_i\doteq\sum_{\ell=1}^s\exp\bigl(\langle\mathsf Qq_i,\mathsf Kq_\ell\rangle/\sqrt d\bigr).
\]
This transformation can be expressed in terms of the empirical probability measure \(
\mu=s^{-1}\sum_{i=1}^s\delta_{q_i}\in \mathcal P(\R^n)\)~--- via the mapping $\mu \mapsto \big[{\id}_{\R^n} + \mathfrak A(\cdot,\mu)\big]_\sharp \mu$, where
\[
\mathfrak A(q,\mu)
\doteq \left(\int_{\R^n}
\exp\left(
\langle \mathsf Qq,\mathsf Kz\rangle/\sqrt d
\right)
\d\mu(z)\right)^{-1}
\int_{\R^n}
\exp\left(
\langle \mathsf Qq,\mathsf Kz\rangle/\sqrt d
\right)
\mathsf Vz\d\mu(z).
\]
For an arbitrary $\mu\in\P$, the formula is meaningful when the denominator is finite and the numerator is Bochner integrable; these conditions are automatic for empirical measures.

A map $\mathcal T\colon\P\to\P$ is called \emph{context-dependent} if, for a suitably regular map $\Gamma\colon Q\times\P\to Q$, it admits the representation
\begin{equation}\label{eq:context-dependent-representation}
\mathcal T(\mu)=\Gamma(\cdot,\mu)_\sharp\mu.
\end{equation}
The class of such maps is described in~\cite{FuruyaDeHoopLassas2025} through their action on empirical measures and a regularity condition w.r.t. variations of the input measure; it contains, in particular, transformations generated by sufficiently regular nonlocal vector fields.

In a Transformer network, the parameterized operation
\(
\mathcal T_{\mathrm u}(\mu)
\doteq
\Gamma_{\mathrm u}(\cdot,\mu)_\sharp\mu
\)
is applied repeatedly in $L$ successive layers, with learnable parameters
$\mathrm u^0,\mathrm u^1,\ldots,\mathrm u^{L-1}$.
The sequence of measures $\mu^k$ therefore satisfies the difference equation
\[
\mu^{k+1}
=
\Gamma_{\mathrm u^k}(\cdot,\mu^k)_\sharp\mu^k,
\quad
k=0,\ldots,L-1; \qquad \mu^0 = \vartheta.
\]
In a residual architecture, one can represent the maps $\Gamma_{\mathrm u}$ in the form $\Gamma_{\mathrm u}(q,\mu)=q+hF(q,\mu,\mathrm u)$, where $h>0$ is the depth step associated with a single layer.
The token vectors and their empirical measures then satisfy the difference equations
\begin{equation}\label{eq:transformer-euler}
q_j^{k+1}
=
q_j^k+h\,F(q_j^k,\mu^k,\mathrm u^k); \quad\quad \mu^{k+1} = \left({\id}_{\R^n} + h \,\! F(\cdot, \mu^k, \mathrm u^k)\right)_\sharp \mu^k.
\end{equation}

Denote by $T=Lh$ the total depth of the network. The first equation in~\eqref{eq:transformer-euler} has the form of an explicit Euler discretization of a controlled ODE, while the sequence $(\mathrm u^k)_{k=0}^{L-1}$ can be seen as a piecewise-constant control
$u_h\colon[0,T]\to U$ given by
\(u_h(t)=\mathrm u^k\) for \(t\in[kh,(k+1)h)\). 
The formal continuous-depth counterpart of~\eqref{eq:transformer-euler}, obtained as $h\downarrow0$ and $L\to\infty$ with $T=Lh$ fixed, is the nonlocal continuity equation~\eqref{eq:CE}; models of this type are studied, e.g., in~\cite{Geshkovski2024,CastinAblinCarrilloPeyre2025,bruno2024emergencemetastableclusteringmeanfield}.

Suppose that the output of the network is described by a map $g\colon\mathcal P(\R^n)\to\R^N$ of the final context $\mu_T$ (or $\mu^L$ in the discrete-depth model).
When the layer dynamics is generated by a finite control-linear family of NVFs as in~\eqref{eq:controlled-nonlocal-field}, steering this observable to prescribed values is an LPC problem, to which Theorem~\ref{thm:transport-projected-rank} applies.

\subsection{Example 1. A bracket-generated direction in a self-attention model}

The preceding discussion identifies a context with an empirical measure whose atoms are the vector representations of the tokens. We now consider the simplest two-token model in $\R^2$, which is a surrogate of the controlled self-attention dynamics in~\cite{GeshkovskiRigolletRuizBalet2026}.

Set $Q=\R^2$ and $\P=\mathcal P(\R^2)$, and let $e_1,e_2$ be the standard basis
of $\R^2$. We take the two token states to be $q_1=e_1$ and $q_2=e_2$, so their
empirical measure is
\(
\mu_0\doteq\frac12\bigl(\delta_{e_1}+\delta_{e_2}\bigr);
\)
we take this measure as the reference distribution $\vartheta$. 

Choose $\chi\in C_c^\infty(\R^2;[0,1])$ such that $\chi=1$ on $B_2(0)$, and put $c_\chi(q)\doteq\chi(q)q$.
Define the self-attention map $\mathfrak A_\chi\colon\R^2\times\P\to\R^2$ by
\[
\mathfrak A_\chi(q,\mu)
\doteq
\left(\int_{\R^2}
\exp\bigl(\langle c_\chi(q),c_\chi(r)\rangle\bigr)
\d\mu(r)\right)^{-1}
\int_{\R^2}
\exp\bigl(\langle c_\chi(q),c_\chi(r)\rangle\bigr)
\,c_\chi(r)\d\mu(r);
\]
when $q$ and the support of $\mu$ lie in $B_2(0)$, it is the usual softmax attention with the identity score matrix, up to the standard scaling factor.

Consider two NVFs $V_1,V_2\colon\R^2\times\P\to\R^2$:
\[
V_1(q,\mu)\doteq e_1,
\qquad
V_2(q,\mu)\doteq M\mathfrak A_\chi(q,\mu), \qquad
M
\doteq
\begin{pmatrix}
0&0\\
1&-1
\end{pmatrix};
\]
$V_1$ is a local VF, while $V_2$ is an NVF depending on the measure through the self-attention map. 

For $U\doteq[-1,1]^2$, define $F\colon\R^2\times\P\times U\to\R^2$ by $F(q,\mu,\mathrm u)\doteq \mathrm u_1V_1(q,\mu)+\mathrm u_2V_2(q,\mu)$, $\mathrm u=(\mathrm u_1,\mathrm u_2)\in U$. As an observable, take $g\colon\P\to\R^2$ defined by $g(\mu)\doteq\langle c_\chi\rangle(\mu)$; the map $c_\chi$ is smooth and bounded together with its derivatives. By Lemma~\ref{lem:smooth-class}, $g$ is smooth, and so are the numerator and the denominator of $\mathfrak A_\chi$, which are NVFs defined by smooth bounded kernels; in particular, their lifts are smooth maps into $C_b(\R^2;\R^2)$ and $C_b(\R^2;\R)$, respectively. Since the denominator is bounded below by $e^{-\|c_\chi\|_\infty^2}>0$ and division by functions bounded away from zero is smooth on $C_b(\R^2;\R)$, the lift of $\mathfrak A_\chi$, and hence that of $V_2$, is smooth. The local VF $V_1$ is constant.
The cutoff is inactive on both atoms of $\vartheta^X$ for $\|X\|_\infty$ sufficiently small. With $k=1$, formula~\eqref{eq:polynomial-calculus} takes the form
\[
(L_1g)(\mu_0)=e_1,
\qquad
(L_2g)(\mu_0)=\frac12 M\bigl(\mathfrak A_\chi(e_1,\mu_0)+\mathfrak A_\chi(e_2,\mu_0)\bigr).
\]
The attention vectors at $e_1$ and $e_2$ are
\[
\mathfrak A_\chi(e_1,\mu_0)=\frac{e\,e_1+e_2}{1+e},
\qquad
\mathfrak A_\chi(e_2,\mu_0)=\frac{e_1+e\,e_2}{1+e}.
\]
Their sum is $e_1+e_2$, and $M(e_1+e_2)=0$. Hence $(L_2g)(\mu_0)=0$, and the first-order projected directions span only the line generated by $e_1$.

For the first commutator, notice that $L_1g(\vartheta^X)=e_1$ whenever $\|X\|_\infty$ is sufficiently small. Consequently, $(L_2L_1g)(\mu_0)=0$.
Let
\[
\mu_\varepsilon
\doteq
({\id}_{\R^2}+\varepsilon e_1)_\sharp\mu_0
=
\frac12\bigl(\delta_{x_1(\varepsilon)}+\delta_{x_2(\varepsilon)}\bigr),
\]
where $x_1(\varepsilon)=(1+\varepsilon)e_1$ and $x_2(\varepsilon)=e_2+\varepsilon e_1$.
By the definition of $L_1$,
\begin{equation}\label{eq:example1-bracket-derivation}
([L_1,L_2]g)(\mu_0)
=
\left.\frac{\d}{\d\varepsilon}\right|_{\varepsilon=0}
M\int_{\R^2}\mathfrak A_\chi(q,\mu_\varepsilon)\d\mu_\varepsilon(q).
\end{equation}
For sufficiently small $|\varepsilon|$, the cutoff is inactive. Denote by $p_i(\varepsilon)$ the softmax weight of $x_1(\varepsilon)$ in the attention vector corresponding to $x_i(\varepsilon)$, $i=1,2$. These weights are
\[
p_1(\varepsilon)
=
\frac{e^{1+\varepsilon}}{1+e^{1+\varepsilon}},
\qquad
p_2(\varepsilon)
=
\frac{e^{\varepsilon-1}}{1+e^{\varepsilon-1}},
\]
with $p_1'(0)=p_2'(0)=\kappa\doteq e/(1+e)^2$. Moreover, the attention vectors satisfy, for $i=1,2$,
\[
\mathfrak A_\chi(x_i(\varepsilon),\mu_\varepsilon)
=x_2(\varepsilon)+p_i(\varepsilon)(e_1-e_2).
\]
The derivative in~\eqref{eq:example1-bracket-derivation} is therefore
\[
([L_1,L_2]g)(\mu_0)
=
M\bigl(e_1+\kappa(e_1-e_2)\bigr)
=
(1+2\kappa)e_2.
\]
The two vectors $(L_1g)(\mu_0)=e_1$ and $([L_1,L_2]g)(\mu_0)=(1+2\kappa)e_2$ span $\R^2$. Hence condition~\eqref{eq:transport-projected-rank} holds with $P_1=L_1$ and $P_2=[L_1,L_2]$, and LPC at $\mu_0$ w.r.t. $g$ follows from Theorem~\ref{thm:transport-projected-rank}.

The calculation of the commutator also shows how the second projected direction arises. At $\mu_0$, the self-attention VF moves the two atoms in opposite directions, so their barycentric contributions cancel. Although the translation generated by $V_1$ leaves the relative position of the atoms unchanged, it affects the attention output in two ways: the direct translation contributes $Me_1=e_2$, while the change in the dot-product attention scores, and hence in the attention weights, contributes $2\kappa e_2$. These two contributions account for the commutator direction computed above.

Since $F$ is bounded, the cutoff remains inactive along these characteristics over sufficiently short control intervals. The LPC conclusion therefore also holds for the uncut dynamics near $\mu_0$.

\section{Symmetry constraints}
\label{sec:equivariance}

We now ask whether dependence of the generating NVFs on the measure is essential for LPC under symmetry restrictions. Consider a system of interacting particles whose velocities are specified by an NVF of the form
\[
V_{\mathrm u}(q,\mu)
=
\int_{\R^n}
K_{\mathrm u}(q-r)\d\mu(r)
\]
with a controlled interaction kernel $K_{\mathrm u}$. If the particles move in a homogeneous space with no distinguished origin, translating the entire configuration should not change their mutual interaction; if the medium is also isotropic, rotating the configuration should simply rotate the resulting velocities in the same reference frame. In order to respect this physical reality, we must admit only those kernels which preserve the mentioned symmetries for every $\mathrm u$.

For example, let $K_{\mathrm u}(z)=\kappa_{\mathrm u}(|z|)z$. The control may change the strength and the distance dependence of the interaction, but every admissible kernel remains translation- and rotation-equivariant (see definitions below). Such symmetry requirements are standard in models of many-body systems. In particular, symmetry-preserving flows have been introduced for physical and chemical particle systems precisely because the underlying probability distributions inherit symmetries from the interaction potential~\cite{KohlerKleinNoe2020}.

Symmetry restrictions modify the very statement of the controllability problem: suppose, for simplicity, that the transported measures have compact supports. If, for example, $K_{\mathrm u}$ is odd, we have
\[
\frac{\d}{\d t}
\int_{\R^n}
q\d\mu_t(q)
=
\int_{\R^n}
\int_{\R^n}
K_{u(t)}(q-r)
\d\mu_t(r)\d\mu_t(q)
=
0,
\]
which shows that the barycenter of $\mu_t$ is preserved by each admissible control. A measure with one barycenter therefore cannot be transported to a measure with another barycenter, no matter how the interaction kernel is varied within this class. Hence full controllability is impossible because the controls are required to remain within the physically admissible family.

A related situation occurs when symmetry is built into the architecture of the model. For three-dimensional point clouds and molecular data, a global rotation or translation changes the coordinate representation without changing the underlying geometric object. Attention mechanisms have therefore been constructed so that rotating or translating the input produces the same rotation or translation of the output for every choice of the trainable parameters. As a concrete instance of such an architecture one may consider SE(3)-Transformers~\cite{FuchsWorrallFischerWelling2020}. The fact that a rotationally symmetric distribution cannot be transformed, within a meaningful class of controls, into an arbitrary nonsymmetric one is well understood (see, e.g.,~\cite{Agrachev2025}); this obstruction is one of the reasons why controllability under symmetry restrictions is studied separately.

\subsection{Equivariance}

Let $G$ be a group acting on $Q$ by smooth diffeomorphisms, with an identity element $e_G$. For $\sigma\in G$, we denote by $\tau_\sigma:Q\to Q$ the corresponding transformation, so that $\tau_{e_G}={\id}_Q$ and $\tau_{\sigma_1\sigma_2}=\tau_{\sigma_1}\circ\tau_{\sigma_2}$. The action on $Q$ induces an action on $\P$ by $\mu\mapsto (\tau_\sigma)_\sharp\mu$.

Recall that a (local) VF $v:Q\to Q$ is called $G$-equivariant if $v(\tau_\sigma(q))=D\tau_\sigma(q)[v(q)]$, $q\in Q$, $\sigma\in G$, which says that the evolution $\dot q_t = v(q_t)$ defined by $v$ commutes with the action of $G$: if $t \mapsto q_t$ is a trajectory starting from $q_0$, then $t\mapsto\tau_\sigma(q_t)$ is the trajectory with initial condition $\tau_\sigma(q_0)$.
Accordingly, an NVF $V:Q\times\P\to Q$ will be called \emph{$G$-equivariant} if
\begin{equation}\label{eq:general-equivariance}
V\bigl(\tau_\sigma(q),(\tau_\sigma)_\sharp\mu\bigr)
=
D\tau_\sigma(q)[V(q,\mu)],
\qquad q\in Q,\quad \mu\in\P,\quad \sigma\in G.
\end{equation}
If an NVF satisfies~\eqref{eq:general-equivariance}, the push-forward under $\tau_\sigma$ of a characteristic solution satisfies the characteristic equation for the transformed initial measure and the same control, provided that the transformed characteristic map has displacement in $C([0,T];\mathscr Q)$. 

Let $\mathcal G\subset\P$ be invariant under the induced action of $G$, i.e., $(\tau_\sigma)_\sharp\mu\in\mathcal G$ for every $\mu\in\mathcal G$ and $\sigma\in G$. A map $g\colon\mathcal G\to\R^N$ is called \emph{$G$-invariant} if $g((\tau_\sigma)_\sharp\mu)=g(\mu)$ for every $\mu\in\mathcal G$ and $\sigma\in G$. 

We are interested in the LPC of such observables when all admissible NVFs are required to satisfy the condition~\eqref{eq:general-equivariance}. We ask whether allowing the NVFs to depend on the current measure changes the available projected directions.

\subsection{Example 2. Translation equivariance}

Let $Q=\R$, and denote by $\P_c\subset\P(\R)$ the set of probability measures with compact support in $\R$. As a single observable we take the variance $g\doteq\mathbb V$ from Section~\ref{ssec:computing-rank}, restricted to $\P_c$. This functional is invariant under the translation $\tau_a(q)\doteq q+a$, i.e., the identity $g((\tau_a)_\sharp\mu)=g(\mu)$ holds for all $\mu\in\P_c$ and $a\in\R$.

An NVF $V\colon\R\times\P_c\to\R$ is called \emph{translation-equivariant} if
\begin{equation}\label{eq:translation-equivariance}
V\bigl(q+a,(\tau_a)_\sharp\mu\bigr)=V(q,\mu),
\qquad \forall \, q,a\in\R,\quad \forall \mu\in\P_c.
\end{equation}
If $V$ is local, i.e. $V(q,\mu)=v(q)$, condition~\eqref{eq:translation-equivariance} reduces to $v(q+a)=v(q)$ for all $q,a\in\R$, so $v$ must be constant. Every trajectory generated by local translation-equivariant VFs is therefore a translation of the initial measure and preserves its variance. Consequently, LPC w.r.t. $g$ does not hold for any family of these VFs.

To see what changes when dependence on the measure is admitted, fix $\vartheta\in\P_c$ and $R>0$ such that $\spt\vartheta\subset[-R,R]$. Choose $\eta\in C_c^\infty(\R)$ equal to $1$ on a neighborhood of $[-2R,2R]$, and define the NVF $V_R\colon\R\times\P_c\to\R$ by
\[
c_\eta(r)\doteq r\eta(r),
\qquad V_R(q,\mu)\doteq c_\eta\bigl(q-\mathbb E(\mu)\bigr).
\]
Since $\mathbb E((\tau_a)_\sharp\mu)=\mathbb E(\mu)+a$, the NVF $V_R$ is translation-equivariant. The identity $\mathbb E(\vartheta^X)=\mathbb E(\vartheta)+\langle X\rangle(\vartheta)$ and the quadratic dependence of $g(\vartheta^X)$ on $X$ show that both functionals are smooth on $\P_c$. The smoothness of $V_R$ follows by composition with $c_\eta$. Since bounded continuous displacements preserve compact support, $\vartheta^X\in\P_c$ for every $X\in\mathscr Q$, as required for the local application of Theorem~\ref{thm:transport-projected-rank}.

For $\mu=\vartheta^X$ with $\|X\|_\infty$ sufficiently small, the cutoff is inactive on $\spt\mu$, and $V_R(q,\mu)=q-\mathbb E(\mu)$ there. Consequently, on the support of $\mu$, the perturbation $q\mapsto q+\varepsilon V_R(q,\mu)=\mathbb E(\mu)+(1+\varepsilon)(q-\mathbb E(\mu))$ is a dilation about the barycenter. Since this transformation multiplies the variance by $(1+\varepsilon)^2$, we obtain
\[
 g\bigl([{\id}_{\R}+\varepsilon V_R(\cdot,\mu)]_\sharp\mu\bigr)
 =(1+\varepsilon)^2g(\mu),
 \qquad (L_{V_R}g)(\mu)=2g(\mu).
\]
By Proposition~\ref{prop:closure}, the projected dynamics is $\dot y=2u(t)y$, $y(0)=g(\vartheta)$. Thus $V_R$ yields LPC at $\vartheta$ w.r.t.\ $g$ if and only if $g(\vartheta)>0$.

\subsection{Example 3. Rotation equivariance}

Identify $Q=\R^2$ with $\C$ and fix $a\in(0,1/2)$. On the unit circle $\mathbb S^1\subset\C$, let $\mu^0$ be the probability measure with angular density $\rho^0(\alpha)\doteq(1+2a\cos\alpha)/(2\pi)$, $\alpha\in[0,2\pi)$. For $\theta\in\R$, we define the rotation $\tau_\theta(q)\doteq e^{i\theta}q$ and the rotated measure $\mu^\theta\doteq(\tau_\theta)_\sharp\mu^0$, whose density is $\rho^\theta(\alpha)\doteq\rho^0(\alpha-\theta)$.

Choose $\zeta\in C_c^\infty([0,\infty))$ such that $\zeta=1$ on $[1/2,3/2]$ and $\zeta=0$ in a neighborhood of the origin. For $s=1,2$, define
\[
z_s(\mu)\doteq\langle h_s\rangle(\mu),
\qquad
h_s(q)\doteq
\begin{cases}
\zeta(|q|)(q/|q|)^s, & q\neq 0,\\
0, & q=0.
\end{cases}
\]
Then $z_1(\mu^\theta)=ae^{i\theta}$, while $z_2(\mu^\theta)=0$.

On the set $\mathcal G\doteq\{\mu\in\P:z_1(\mu)\neq 0\}$, consider the observable
\[
g(\mu)\doteq
\left(
\mathrm{Re}\frac{z_2(\mu)}{z_1(\mu)^2},
\mathrm{Im}\frac{z_2(\mu)}{z_1(\mu)^2}
\right)\in\R^2.
\]
By the identities $z_s((\tau_\theta)_\sharp\mu)=e^{is\theta}z_s(\mu)$, $s=1,2$, the map $g$ is rotation invariant; in particular, $g(\mu^\theta)=0_{\R^2}$. In this case, condition~\eqref{eq:general-equivariance} takes the form
\[
V\bigl(e^{i\theta}q,(\tau_\theta)_\sharp\mu\bigr)
=
e^{i\theta}V(q,\mu),
\qquad q\in\C,\quad \mu\in\P,\quad \theta\in\R.
\]

If $V$ is local, i.e. $V(q,\mu)=v(q)$, then condition~\eqref{eq:general-equivariance} reduces to $v(e^{i\theta}q)=e^{i\theta}v(q)$. Consider a control-linear system generated by finitely many smooth bounded $SO(2)$-equivariant VFs $v_j$. For $q\neq0$, every such VF has the form $v_j(q)=\beta_j(|q|)q$ with a complex function $\beta_j$.

Fix $T>0$ and a control $u\in\mathcal U_T$, and put $\beta_u(t,r)\doteq\sum_{j=1}^m u_j(t)\beta_j(r)$. All characteristics starting from $\mathbb S^1$ have the form $t\mapsto\lambda_t q$, where
\[
\dot\lambda_t=\beta_u(t,|\lambda_t|)\lambda_t,
\qquad \lambda_0=1.
\]
The solution satisfies $\lambda_t=\exp\!\left(\int_0^t\beta_u(s,|\lambda_s|)\d s\right)\neq0$, and, for $s=1,2$, the corresponding moments are
\[
 z_s(\mu_t^u)
 =\zeta(|\lambda_t|)\left(\frac{\lambda_t}{|\lambda_t|}\right)^s z_s(\mu^\theta).
\]
Since $z_2(\mu^\theta)=0$, the displayed formula gives $z_2(\mu_t^u)=0$ along every admissible trajectory. Thus $g(\mu_T^u)=0_{\R^2}$ for any terminal time $T>0$ and every admissible control for which $\mu_T^u\in\mathcal G$, so no system generated by local $SO(2)$-equivariant VFs is LPC at $\mu^\theta$ w.r.t. $g$.

The situation changes when the VFs are allowed to depend on the current measure, subject to the same equivariance condition. Choose $\eta\in C_c^\infty([0,\infty))$ such that $\eta=1$ in a neighborhood of $1$ and $\eta=0$ in a neighborhood of the origin, and define
\[
\Omega(q)\doteq
\begin{cases}
\eta(|q|)\dfrac{q^2}{|q|^2}, & q\neq 0,\\
0, & q=0,
\end{cases}
\qquad
\Psi(q)\doteq\eta(|q|)\,iq.
\]
Choose also $\chi\in C_c^\infty([0,\infty))$, with $0\leq\chi\leq1$, such that $\chi=0$ in a neighborhood of the origin and $\chi=1$ in a neighborhood of $a$, and put
\[
\Lambda(\mu)\doteq
\begin{cases}
\displaystyle
\chi(|z_1(\mu)|)\frac{z_1(\mu)^2}{|z_1(\mu)|^2},
& z_1(\mu)\neq 0,\\[0.8em]
0, & z_1(\mu)=0.
\end{cases}
\]
Then $\Omega(e^{i\theta}q)=e^{2i\theta}\Omega(q)$, $\Psi(e^{i\theta}q)=e^{i\theta}\Psi(q)$, and $\Lambda((\tau_\theta)_\sharp\mu)=e^{2i\theta}\Lambda(\mu)$.

Consider the NVFs $V_1,V_2\colon\C\times\P\to\C$ defined by
\[
V_1(q,\mu)\doteq
\mathrm{Re}\bigl(\overline{\Lambda(\mu)}\,\Omega(q)\bigr)\Psi(q),
\qquad
V_2(q,\mu)\doteq
\mathrm{Im}\bigl(\overline{\Lambda(\mu)}\,\Omega(q)\bigr)\Psi(q).
\]
They satisfy the same $SO(2)$-equivariance condition. The cutoff functions remove the singularities at $q=0$ and $z_1(\mu)=0$. In particular, $\Omega$, $\Psi$, $h_s$, and the map $z\mapsto\chi(|z|)z^2/|z|^2$, extended by zero at $z=0$, are smooth with bounded derivatives. It follows from Lemma~\ref{lem:smooth-class} that $z_1$, $z_2$, $V_1$ and $V_2$ are smooth, while $g$ is smooth on $\mathcal G$. For $U=[-1,1]^2$, we define the controlled NVF $F\colon\C\times\P\times U\to\C$ by $F(q,\mu,\mathrm u)\doteq\mathrm u_1V_1(q,\mu)+\mathrm u_2V_2(q,\mu)$.

Fix $\theta_0\in\R$. Since $\Lambda(\mu^{\theta_0})=e^{2i\theta_0}$, on $\mathbb S^1$ we have
\[
V_1(e^{i\alpha},\mu^{\theta_0})
=
\cos\bigl(2(\alpha-\theta_0)\bigr)ie^{i\alpha},
\qquad
V_2(e^{i\alpha},\mu^{\theta_0})
=
\sin\bigl(2(\alpha-\theta_0)\bigr)ie^{i\alpha}.
\]

We take $\mu^{\theta_0}\in\mathcal G$ as the reference measure and write $L_j\doteq L_{V_j}$, $j=1,2$. At $\mu^{\theta_0}$, where $z_2(\mu^{\theta_0})=0$, we apply the quotient rule to $g$ and obtain
\begin{equation}\label{eq:example3-Ljg}
(L_jg)(\mu^{\theta_0})
=
\left(
\mathrm{Re}\left[(ae^{i\theta_0})^{-2}(L_jz_2)(\mu^{\theta_0})\right],
\mathrm{Im}\left[(ae^{i\theta_0})^{-2}(L_jz_2)(\mu^{\theta_0})\right]
\right).
\end{equation}
By Proposition~\ref{prop:polynomial-calculus},
\begin{equation}\label{eq:example3-Ljz2}
(L_jz_2)(\mu^{\theta_0})
=
\int_\C Dh_2(q)
\bigl[V_j(q,\mu^{\theta_0})\bigr]
\d\mu^{\theta_0}(q).
\end{equation}
For $q=e^{i\alpha}\in\mathbb S^1$ and $\omega\in\R$, the curve $t\mapsto e^{i(\alpha+\omega t)}$ remains on $\mathbb S^1$, where $\zeta=1$. Consequently,
\[
Dh_2(e^{i\alpha})[\omega i e^{i\alpha}]=2i\omega e^{2i\alpha}.
\]
Let $\omega_1(\alpha)\doteq\cos\bigl(2(\alpha-\theta_0)\bigr)$ and $\omega_2(\alpha)\doteq\sin\bigl(2(\alpha-\theta_0)\bigr)$. Using the expressions for $V_1$ and $V_2$ on $\mathbb S^1$, we compute the integral in~\eqref{eq:example3-Ljz2}:
\[
(L_jz_2)(\mu^{\theta_0})
=2i\int_0^{2\pi}e^{2i\alpha}\omega_j(\alpha)\rho^{\theta_0}(\alpha)\d\alpha
=\begin{cases}
 ie^{2i\theta_0},&j=1,\\
 -e^{2i\theta_0},&j=2.
\end{cases}
\]
Returning to~\eqref{eq:example3-Ljg}, the two projected directions are
\[
(L_1g)(\mu^{\theta_0})
=
\left(0,\frac1{a^2}\right),
\qquad
(L_2g)(\mu^{\theta_0})
=
\left(-\frac1{a^2},0\right).
\]
These vectors span $\R^2$, so condition~\eqref{eq:transport-projected-rank} holds already with $P_1=L_1$ and $P_2=L_2$. By Theorem~\ref{thm:transport-projected-rank}, the system generated by $V_1$ and $V_2$ is LPC at $\mu^{\theta_0}$ w.r.t. $g$.

Thus, subject to the same rotation-equivariance constraint, the observable $g$ is fixed by every admissible local system, whereas the full projected rank is attained by the two NVFs above without involving higher commutators.

\subsection{Context-independent maps}

The dependence of the transformation on the current measure is commonly viewed as an important source of the ``expressive power'' of Transformers. At the same time, the numerical price of nonlocality is typically high: when the transformation in~\eqref{eq:context-dependent-representation} is generated by an NVF, the characteristic trajectories are generally coupled through the evolving measure and must be computed simultaneously. For a local VF, by contrast, the continuity equation is linear and the characteristics can be tracked independently.

It is therefore natural to ask when a map $\mathcal T\colon\P\to\P$ admits the form $\mathcal T(\mu)=S_\sharp\mu$, $\mu\in\P$, with a continuous map $S\colon Q\to Q$ independent of the measure. We call such a transformation \emph{context-independent}. The following proposition characterizes these transformations among maps on $\P$ continuous in the narrow topology.

\begin{proposition}\label{prop:context-independent-map}
Let \(\mathcal T\colon\P\to\P\) be narrowly continuous. Then the following
assertions are equivalent.
\begin{enumerate}
\item There exists a map \(S\in C(Q;Q)\) such that
\(\mathcal T(\mu)=S_\sharp\mu\) for every \(\mu\in\P\).

\item The map \(\mathcal T\) sends Dirac measures to Dirac measures, that is, for every $q\in Q$, there exists a unique point $\tilde q\in Q$ such that
\(
\mathcal T(\delta_q)=\delta_{\tilde q}.
\)
Moreover $\mathcal T$ is affine, i.e.,
\[
\mathcal T\bigl((1-\lambda)\mu+\lambda\nu\bigr)
=
(1-\lambda)\mathcal T(\mu)+\lambda\mathcal T(\nu)
\]
for all \(\mu,\nu\in\P\) and \(\lambda\in[0,1]\).
\end{enumerate}
\end{proposition}

Finite compositions preserve the class of transformations characterized by Proposition~\ref{prop:context-independent-map}; it also contains the flow maps of local VFs whenever their flows are defined on all of $Q$. For the action introduced above, we call a map $\mathcal T\colon\P\to\P$ $G$-equivariant if $\mathcal T((\tau_\sigma)_\sharp\mu)=(\tau_\sigma)_\sharp\mathcal T(\mu)$ for every $\mu\in\P$ and $\sigma\in G$. If $\mathcal T=S_\sharp$ is $G$-equivariant, then $S$ is $G$-equivariant, as follows by evaluating this identity on Dirac measures.

This observation strengthens the two examples above. Under translations on $\R$, an equivariant map has the form $S(q)=q+b$, so the variance is preserved. Under rotations on $\C$, an equivariant map satisfies $S(e^{i\alpha})=e^{i\alpha}S(1)$ on the unit circle; therefore the second harmonic of $S_\sharp\mu^{\theta}$ remains zero. Whenever the transformed measure belongs to $\mathcal G$, its value of $g$ is consequently $0_{\R^2}$. Thus, in both examples, the obstruction is not restricted to flows of local VFs but extends to all continuous context-independent equivariant transformations.

\section{Conclusion}

This work suggests a ``practical'' reading of the controllability problem for nonlocal continuity equations on the space of probability measures. The key technical ingredient is the representation of nonlocal characteristics (and the respective curves of measures) by solutions of an ODE on a fixed Banach space; this interpretation is of independent interest for the analysis on the space of measures and mean-field control, beyond the controllability theory.

As with most results of their kind, Theorems~\ref{theorem} and \ref{thm:transport-projected-rank} give neither a useful estimate of the reachable neighborhood and the time horizon $T$, nor an a priori bound on the depth of Lie brackets required to obtain the full projected rank. A quantitative version of LPC, relating the size of the reachable neighborhood and the control time to the projected Lie algebra, remains an open problem (cf. Remark~\ref{re:numerical}). Notice also that the rank condition obtained here is only sufficient. It would be interesting to understand to what extent it can be complemented by necessary conditions.
 
Although our examples indicate that the dependence of an NVF on the current measure may be essential, it is not clear how general this effect is. A natural question is to characterize, for a given group action and a given invariant observation, when the dependence on the measure enlarges the projected Lie algebra.

Finally, the connection with Transformer architectures remains qualitative. For applications, one should understand how LPC or related properties survive a finite-depth discretization and how the corresponding bracket mechanisms are reflected in finite compositions of self-attention layers.

\bmhead{Acknowledgements}

M.~Staritsyn acknowledges financial support of the Research Center for Systems and Technologies (UID/00147) and the Associate Laboratory ARISE (LA/\allowbreak P/\allowbreak 0112/\allowbreak 2020, DOI:~10.\allowbreak 54499/\allowbreak LA/\allowbreak P/\allowbreak 0112/\allowbreak 2020) funded by Funda\c{c}\~ao para a Ci\^encia e a Tecnologia, I.P./MECI through national funds and the grant 2023.\allowbreak 05910.\allowbreak CEECIND/\allowbreak CP2834/\allowbreak CT0017, DOI:~10.\allowbreak 54499/\allowbreak 2023.\allowbreak 05910.\allowbreak CEECIND/\allowbreak CP2834/\allowbreak CT0017.

\appendix

\section{Banach manifolds, vector fields and Lie brackets}\label{app:manifold-notation}

We recall the basic constructions for a smooth Banach manifold modeled on a real Banach space $E$; see~\cite[pp.~22--154]{Lang1999}. A chart $(\mathcal O,\psi)$ identifies $\mathcal O\subset\mathcal M$ with an open subset of $E$, and its differential $D\psi(x)\colon\mathrm T_x\mathcal M\to E$, $x\in\mathcal O$, is a linear isomorphism. The coordinate representation of $v\in\mathrm T_x\mathcal M$ is $v_\psi\doteq D\psi(x)[v]$.

A VF on $\mathcal M$ is a map $F\colon\mathcal M\to\mathrm T\mathcal M$ with $F(x)\in\mathrm T_x\mathcal M$ for every $x\in\mathcal M$. In a chart $(\mathcal O,\psi)$ its coordinate representation is $F^\psi(\xi)\doteq D\psi(x)[F(x)]$, where $x=\psi^{-1}(\xi)$.
The derivative of a smooth function $\phi$ on $\mathcal M$ along $F$ is
\[
(F\phi)(x)
=
D(\phi\circ\psi^{-1})(\xi)[F^\psi(\xi)],
\qquad \xi=\psi(x),
\]
and is independent of the chosen chart.

If $F_1$ and $F_2$ are smooth VFs, their \emph{Lie bracket} $[F_1,F_2]_{\mathcal M}$ is determined by
\[
D\psi(x)\bigl[[F_1,F_2]_{\mathcal M}(x)\bigr]
=
DF_2^\psi(\xi)[F_1^\psi(\xi)]
-
DF_1^\psi(\xi)[F_2^\psi(\xi)],
\qquad \xi=\psi(x).
\]
The bracket is independent of the chart~\cite[Proposition~1.3, pp.~118--119]{Lang1999}. Our sign convention is $[F_1,F_2]_{\mathcal M}\phi=F_1(F_2\phi)-F_2(F_1\phi)$, as stated in Section~\ref{sec:prelim}.

\section{Proofs}

\subsection{Proof of Theorem~\ref{theorem}}

Let $e^{t v_i}$ denote the local flow of $v_i$, $i=1,\ldots,m$. Given local diffeomorphisms $\Theta_1$ and $\Theta_2$, put $[\Theta_1,\Theta_2]\doteq \Theta_2^{-1}\circ\Theta_1^{-1}\circ\Theta_2\circ\Theta_1$. If $B$ is a bracket monomial, let $k_B$ be the number of occurrences of the generators $v_i$ in $B$, and define the commutator curves recursively by
\[
\mathcal C_{v_i}(\tau)=e^{\tau v_i},
\qquad
\mathcal C_{[B_1,B_2]}(\tau)
=[\mathcal C_{B_1}(\tau),\mathcal C_{B_2}(\tau)].
\]
By the commutator theorem for smooth VFs on Banach manifolds~\cite[Theorem~5.1]{KipkaLedyaev2015} (see also~\cite{Mauhart1992CommutatorsOF}), in a chart $(\mathcal O,\psi)$ one has
\[
\begin{aligned}
(\psi\circ\mathcal C_B(\tau)\circ\psi^{-1})(\xi)
&=
\xi+\tau^{k_B}B^\psi(\xi)+\tau^{k_B+1}R_B(\tau,\xi),\\
\frac{\partial}{\partial\tau}
(\psi\circ\mathcal C_B(\tau)\circ\psi^{-1})(\xi)
&=
k_B\tau^{k_B-1}B^\psi(\xi)+\tau^{k_B}\widetilde R_B(\tau,\xi),
\end{aligned}
\]
where the remainders, together with their first derivatives w.r.t.~$\xi$, are continuous and locally bounded. The inverse commutator curve has the same expansions with the opposite leading terms.

When $B=v_i$, put $\Phi_B(t)=e^{t v_i}$. If $k_B>1$, put
\[
\Phi_B(t)
\doteq
\begin{cases}
\mathcal C_B(t^{1/k_B}), & t\geq0,\\[0.3em]
\mathcal C_B((-t)^{1/k_B})^{-1}, & t<0.
\end{cases}
\]
Substituting $\tau=|t|^{1/k_B}$ in the preceding expansions, we obtain
\[
\begin{aligned}
\psi\circ\Phi_B(t)\circ\psi^{-1}(\xi)
&=
\xi+tB^\psi(\xi)+O(|t|^{1+1/k_B}),\\
\partial_t(\psi\circ\Phi_B(t)\circ\psi^{-1})(\xi)
&=
B^\psi(\xi)+O(|t|^{1/k_B}),
\end{aligned}
\]
locally uniformly together with the first derivative w.r.t.~$\xi$. Hence $\Phi_B$ is continuously differentiable near $(0,x_0)$, $\Phi_B(0)(x)=x$, and $\partial_t\Phi_B(0)(x)=B(x)$.

Let $L_B$ denote the number of elementary flow segments in the commutator representation of $B$. Each factor of $\mathcal C_B$ or $\mathcal C_B^{-1}$ is a flow of some $\pm v_i$ and therefore corresponds to an admissible constant control, since $[-1,1]^m\subset U$. Thus $\Phi_B(t)$ is realized by a piecewise-constant control of total duration $L_B|t|^{1/k_B}$.

Each $W_i$ in~\eqref{c:pr} is a finite linear combination of bracket monomials. Hence one can choose monomials $B_1,\ldots,B_N$ such that the vectors $Dg(x_0)[B_j(x_0)]$, $j=1,\ldots,N$, form a basis of $\R^N$. Put $\Phi_j\doteq\Phi_{B_j}$ and, on a sufficiently small neighborhood $\mathcal W_0$ of $0_{\R^N}$, define
\[
\bm\Phi(t_1,\ldots,t_N)
\doteq
\Phi_N(t_N)\circ\cdots\circ\Phi_1(t_1)(x_0).
\]
The map $\bm\Phi$ is continuously differentiable and $\partial_{t_j}\bm\Phi(0_{\R^N})=B_j(x_0)$. Therefore, for $\mathcal F\doteq g\circ\bm\Phi$,
\[
\partial_{t_j}\mathcal F(0_{\R^N})
=
Dg(x_0)[B_j(x_0)],
\qquad j=1,\ldots,N,
\]
and $D\mathcal F(0_{\R^N})$ is invertible.

Choose $\delta>0$ so small that all concatenations with $|t_j|<\delta$ stay in the domain of local well-posedness and the total control time satisfies
\(
 T\doteq\sum_{j=1}^N L_{B_j}\delta^{1/k_{B_j}}<T_0.
\)
After decreasing the neighborhood if necessary, the inverse function theorem shows that $\mathcal F$ maps a neighborhood of $0_{\R^N}$ onto a neighborhood $\mathcal N$ of $g(x_0)$. For every $y\in\mathcal N$, the corresponding point $\bm\Phi(t_1,\ldots,t_N)$ is reached by a finite concatenation of admissible constant controls of total duration at most $T$. If this duration is smaller than $T$, we extend the control by zero on the remaining interval. Hence $g(x^u(T;x_0))=y$ for some piecewise-constant $u\in\mathcal U_T$. The control horizon $T$ tends to zero with $\delta$.

\subsection{Proof of Lemma~\ref{lem:smooth-class}}

Given $r=(r_1,\ldots,r_k)\in Q^k$ and $X,H\in\mathscr Q$, put
\[
\bm r^X\doteq
(\widetilde X(r_1),\ldots,\widetilde X(r_k)),
\qquad
\bm H(r)\doteq(H(r_1),\ldots,H(r_k)).
\]
Consider the map $f\colon\P\to E$ defined by $f(\mu)=\langle\phi\rangle(\mu)$ for a smooth kernel $\phi\colon Q^k\to E$ bounded together with all its derivatives. The lift $\widehat f\colon\mathscr Q\to E$ is given by
\[
\widehat f(X)=\int_{Q^k}\phi(\bm r^X)\d\vartheta^{\otimes k}(r).
\]
Its derivatives of order $j\geq1$ are
\[
D^j\widehat f(X)[H_1,\ldots,H_j]
=\int_{Q^k}D^j\phi(\bm r^X)
[\bm H_1(r),\ldots,\bm H_j(r)]\d\vartheta^{\otimes k}(r).
\]
Indeed, the first-order Taylor remainder is uniformly bounded in $r$ by $\frac12\|D^2\phi\|_\infty\|H\|_\infty^2$.
To differentiate the displayed expression once more, we apply the same estimate to $D^j\phi$, with $D^{j+2}\phi$ in place of $D^2\phi$. Moreover,
\[
\bigl\|D^j\widehat f(X)[H_1,\ldots,H_j]\bigr\|
\leq\|D^j\phi\|_\infty\prod_{i=1}^j\|H_i\|_\infty,
\]
where $Q^k$ is equipped with the maximum product norm. Since $D^{j+1}\phi$ is bounded, $D^j\phi$ is uniformly continuous, and the above formula shows that $D^j\widehat f$ is continuous in the operator norm. Thus $\widehat f$ is smooth on $\mathscr Q$.

To differentiate the integral defining the lifted NVF, set $\bm r^X_q\doteq(\widetilde X(q),\bm r^X)$ and $\bm H_q(r)\doteq(H(q),\bm H(r))$. Differentiating under the integral sign, we have:
\[
\begin{aligned}
&D^j\!\left[\,X\mapsto\int_{Q^k}\phi(\bm r^X_{\cdot})\d\vartheta^{\otimes k}(r)\right]
(X)[H_1,\ldots,H_j](q)\\
&=\int_{Q^k}D^j\phi(\bm r^X_q)
[\bm H_{1,q}(r),\ldots,\bm H_{j,q}(r)]\d\vartheta^{\otimes k}(r) \qquad \forall j\geq 1.
\end{aligned}
\]
These estimates are uniform in $q\in Q$. Consequently, the derivatives are bounded multilinear maps into $C_b(Q;E)$ and are continuous in the corresponding operator norm. The standard differentiation rules give the assertions on sums and compositions, as well as those on products and quotients of scalar-valued functionals on domains where the denominator stays away from zero.

\subsection{Proof of Proposition~\ref{prop:local-lift}}

We take a common neighborhood $\mathcal D$ on which all $\widehat V_j$ are defined and smooth. By continuity of $D\widehat V_j$ at $0_{\mathscr Q}$, one may choose $r>0$ and $M>0$ such that $\overline B_r(0_{\mathscr Q})\subset\mathcal D$ and, for all $X,Y$ in this ball,
\[
\|\widehat V_j(X)\|_\infty\leq M,
\qquad
\|\widehat V_j(X)-\widehat V_j(Y)\|_\infty
\leq M\|X-Y\|_\infty.
\]
Put $C\doteq\sup_{\mathrm u\in U,\,1\leq j\leq m}|\mathrm u_j|$. Along any continuous curve $t\mapsto X_t$ in the ball, the maps $t\mapsto\widehat V_j(X_t)$ are continuous in $\mathscr Q$. Hence $t\mapsto\sum_{j=1}^m u_j(t)\widehat V_j(X_t)$ is strongly measurable and Bochner integrable. The right-hand side of~\eqref{eq:lifted-control-system} is Lipschitz in $X$ uniformly with respect to $u$ and its norm does not exceed $mCM$ on the ball. Choose $T_0>0$ so that $mCMT_0<r$. The existence and uniqueness theorem for Carath\'eodory ODEs in Banach spaces~\cite[pp.~4--17]{Deimling1977} then ensures that, for every $u\in\mathcal U_T$, $T\leq T_0$, there is a unique absolutely continuous $\mathscr Q$-valued solution on $[0,T]$. In particular, $\|X_t^u\|_\infty\leq mCMt<r$, independently of the admissible control.

Evaluating the integral form of~\eqref{eq:lifted-control-system} at $q\in Q$ and using the definition of $\widehat V_j$, we get
\[
X_t^u(q)=\int_0^t F\bigl(q+X_s^u(q),\vartheta^{X_s^u},u(s)\bigr)\d s.
\]
By~\eqref{eq:lift-flow-reconstruction}, $\vartheta^{X_s^u}=(\Xi_s^{u,\vartheta})_\sharp\vartheta$; consequently, the characteristic integral equation~\eqref{eq:nlf} holds. Conversely, suppose that a map satisfying~\eqref{eq:nlf} has displacement $t\mapsto\Xi_t^{u,\vartheta}-{\id}_Q$ in $C([0,T];\mathscr Q)$. Up to its first exit from $B_r(0_{\mathscr Q})$, this displacement solves the same $\mathscr Q$-valued integral equation as $X^u$. Local uniqueness and the estimate $\|X_t^u\|_\infty<r$ exclude such an exit before $T$.

It remains to verify the distributional identity. Let $\phi$ be a test function in~\eqref{eq:weak-formulation}. Using the definition of push-forward, we express the increment over $[s,t]\subset[0,T]$ as
\[
\langle\phi\rangle(\mu_t^u)-\langle\phi\rangle(\mu_s^u)
=\int_Q\bigl(\phi(\Xi_t^u(q))-\phi(\Xi_s^u(q))\bigr)\d\vartheta(q).
\]
Integrating the chain rule along the absolutely continuous characteristic curve issued from $q\in Q$, we find
\[
\phi(\Xi_t^u(q))-\phi(\Xi_s^u(q))
=\int_s^t D\phi(\Xi_\tau^u(q))
\bigl[F(\Xi_\tau^u(q),\mu_\tau^u,u(\tau))\bigr]\d\tau.
\]
Since $D\phi$ is bounded and the $\widehat V_j$ are uniformly bounded along $X^u$, the integrand is bounded by $mCM\|D\phi\|_\infty$. Fubini's theorem and~\eqref{eflowsol} now give
\[
\begin{aligned}
\langle\phi\rangle(\mu_t^u)-\langle\phi\rangle(\mu_s^u)=\int_s^t\int_Q D\phi(q)
\bigl[F(q,\mu_\tau^u,u(\tau))\bigr]\d\mu_\tau^u(q)\d\tau.
\end{aligned}
\]
Thus, $t\mapsto\langle\phi\rangle(\mu_t^u)$ is absolutely continuous, and its derivative is given by~\eqref{eq:weak-formulation} for a.e.~$t$.

\subsection{Proof of Lemma~\ref{lem:LV-chain-rule}}

By~\eqref{eq:lift-flow-reconstruction}, $f(\mu_t^u)=\widehat{f}(X_t^u)$. Since $X^u$ is Lipschitz and solves~\eqref{eq:lifted-control-system} for a.e.~$t$, the chain rule in $\mathscr Q$, together with~\eqref{eq:derivation-lift}, yields
\[
\frac{\d}{\d t}f(\mu_t^u)
=
D\widehat{f}(X_t^u)\!\left[\sum_{j=1}^m u_j(t)\widehat V_j(X_t^u)\right]
=
\sum_{j=1}^m u_j(t)(L_jf)(\mu_t^u).
\]

\subsection{Proof of Proposition~\ref{prop:bracket}}

Fix $X\in\mathcal D$ and put $\mu\doteq\vartheta^X$. Since $[{\id}_Q+\varepsilon V(\cdot,\mu)]\circ\widetilde X=\widetilde X+\varepsilon\widehat V(X)$, we have $[{\id}_Q+\varepsilon V(\cdot,\mu)]_\sharp\mu=\vartheta^{X+\varepsilon\widehat V(X)}$. Hence, for every $q\in Q$,
\[
\widehat W\bigl(X+\varepsilon\widehat V(X)\bigr)(q)
=W\Bigl(\widetilde X(q)+\varepsilon V\bigl(\widetilde X(q),\mu\bigr),\bigl[{\id}_Q+\varepsilon V(\cdot,\mu)\bigr]_\sharp\mu\Bigr).
\]
The left-hand side is differentiable in $\varepsilon$, since $\widehat W$ is smooth and evaluation at $q$ is a continuous linear map on $\mathscr Q$; its derivative at $\varepsilon=0$ is $D\widehat W(X)[\widehat V(X)](q)$. The right-hand side depends on $q$ only through $p=\widetilde X(q)$. Therefore $\nabla_VW(p,\mu)$ exists for every $p$ in the range of $\widetilde X$, its value does not depend on the choice of $q$ with $\widetilde X(q)=p$, and
\[
\nabla_VW\bigl(\widetilde X(q),\mu\bigr)=D\widehat W(X)[\widehat V(X)](q).
\]
Note that every $p\in\R^n$ lies in the range of $\widetilde X$. Indeed, the continuous map $q\mapsto p-X(q)$ sends $\R^n$, and in particular the closed ball of radius $\|X\|_\infty$ around $p$, into this ball; by Brouwer's theorem it has a fixed point $q$, i.e., $\widetilde X(q)=p$. Interchanging $V$ and $W$ and subtracting, we obtain
\[
\bigl(\nabla_VW-\nabla_WV\bigr)\bigl(\widetilde X(q),\vartheta^X\bigr)=\bigl(D\widehat W(X)[\widehat V(X)]-D\widehat V(X)[\widehat W(X)]\bigr)(q)=[\widehat V,\widehat W](X)(q).
\]

\subsection{Proof of Theorem~\ref{thm:transport-projected-rank}}

Choose $P_1,\ldots,P_N$ satisfying the rank condition~\eqref{eq:transport-projected-rank}, and let $\widehat P_i$ be the corresponding linear combinations of iterated commutators of the lifted VFs $\widehat V_j$. Identity~\eqref{eq:Lie-derivation-lift} gives the equality
\[
(P_i g)(\vartheta)
=
D\widehat g(0_{\mathscr Q})[\widehat P_i(0_{\mathscr Q})],
\qquad i=1,\ldots,N.
\]
Thus~\eqref{eq:transport-projected-rank} is exactly the projected rank condition~\eqref{c:pr} for the lifted system~\eqref{eq:lifted-control-system} at $0_{\mathscr Q}$ with the observation $\widehat g$.

Take a common neighborhood $\mathcal D$ of $0_{\mathscr Q}$ on which $\widehat g$ and all $\widehat V_j$ are defined and smooth. By Proposition~\ref{prop:local-lift}, the lifted trajectories remain in $\mathcal D$ up to some $T_0>0$, uniformly in the control. Theorem~\ref{theorem} then implies that, for every sufficiently small $T>0$, the set of terminal values $\widehat g(X_T^u)$ obtained by piecewise-constant controls contains a neighborhood of $\widehat g(0_{\mathscr Q})=g(\vartheta)$.

For the same controls, we have $\mu_T^u=\vartheta^{X_T^u}$ with $X_T^u\in\mathcal D$ by Proposition~\ref{prop:local-lift}; hence $g(\mu_T^u)$ is defined and equals $\widehat g(X_T^u)$. Therefore the same neighborhood of $g(\vartheta)$ is reached by characteristic solutions of~\eqref{eq:CE}, which proves LPC.

\subsection{Proof of Proposition~\ref{prop:closure}}

For sufficiently small $T>0$, Proposition~\ref{prop:local-lift}, applied on a common neighborhood where $\widehat f$ and all $\widehat V_j$ are smooth, gives $\mu_t^u=\vartheta^{X_t^u}$ with $X_t^u\in\mathcal D$ for $t\in[0,T]$. Hence~\eqref{eq:closure-relation} holds at $\mu_t^u$. Applying Lemma~\ref{lem:LV-chain-rule}, we find, for a.e.~$t$,
\[
\frac{\d}{\d t}f(\mu_t^u)
=
\sum_{j=1}^m
u_j(t)(L_jf)(\mu_t^u)
=
\sum_{j=1}^m
u_j(t)a_j(f(\mu_t^u)).
\]
Thus $y(t)=f(\mu_t^u)$ solves~\eqref{eq:closed-system}. Uniqueness follows since the maps $a_j$ are smooth, hence locally Lipschitz; the equality of the reachable values of $g=h\circ f$ and $h(y)$ follows immediately.

\subsection{Proof of Proposition~\ref{prop:polynomial-calculus}}

Fix $\mu$ and put $\mu_\varepsilon\doteq[{\id}_Q+\varepsilon V(\cdot,\mu)]_\sharp\mu$. Then
\[
\langle\phi\rangle(\mu_\varepsilon)
=\int_{Q^k}
\phi\bigl(q_1+\varepsilon V(q_1,\mu),\ldots,
q_k+\varepsilon V(q_k,\mu)\bigr)\d\mu^{\otimes k}.
\]
The difference quotients are bounded by $k\|D\phi\|_\infty\|V(\cdot,\mu)\|_\infty$. We may therefore differentiate under the integral sign by dominated convergence and obtain~\eqref{eq:polynomial-calculus}.

Take $V(q,\mu)=\langle K(q,\cdot)\rangle(\mu)$. Since the integrands are bounded, Fubini's theorem applied to~\eqref{eq:polynomial-calculus} gives $(L_V\langle\phi\rangle)(\mu)=\langle\partial_K\phi\rangle(\mu)$.
Since $\partial_K\phi$ has $k+1$ arguments, the same calculation can be repeated with $V_i$ or $V_j$. Subtracting the two compositions, we obtain~\eqref{eq:interaction-commutator}. Each iteration adds one argument to the kernel. For smooth functions of finitely many integral terms, the assertion follows by the chain and Leibniz rules.

\subsection{Proof of Proposition~\ref{prop:context-independent-map}}

The first assertion implies the second one. Conversely, suppose that $\mathcal T$ satisfies the second assertion. For every $q\in Q$, there is then a unique point $\tilde q\in Q$ for which
\(
\mathcal T(\delta_q)=\delta_{\tilde q},
\)
and we set $S(q)\doteq\tilde q$. To prove that $S$ is continuous, let $q_k\to q$ in $Q$. Then $\delta_{q_k}\to\delta_q$ narrowly, and by the continuity of $\mathcal T$, $\delta_{S(q_k)}\to\delta_{S(q)}$ narrowly. Testing with the bounded continuous function $r\mapsto\min\{1,\|r-S(q)\|\}$, we get $S(q_k)\to S(q)$. Thus $S\in C(Q;Q)$.

Let \(\mu=\sum_{j=1}^s\lambda_j\delta_{q_j}\) be a finitely supported measure. By the affinity of \(\mathcal T\),
\[
\mathcal T(\mu)
=
\sum_{j=1}^s\lambda_j\mathcal T(\delta_{q_j})
=
\sum_{j=1}^s\lambda_j\delta_{S(q_j)}
=
S_\sharp\mu.
\]
Since $Q$ is separable, finitely supported probability measures are narrowly dense in $\P$, and the narrow topology is metrizable. Hence, for an arbitrary $\mu\in\P$, there are finitely supported $\mu_k\to\mu$ narrowly. Since $S$ is continuous, $S_\sharp\mu_k\to S_\sharp\mu$ narrowly by~the continuous mapping theorem; see, e.g.,~\cite[Chapter~11]{dudley:89:real}. Passing to the limit in $\mathcal T(\mu_k)=S_\sharp\mu_k$ and using the continuity of $\mathcal T$, we obtain $\mathcal T(\mu)=S_\sharp\mu$ for every $\mu\in\P$.

\bibliography{Lie-NLCL-MCSS_Ready_for_submission}

@article{AgrachevSarychev2006,
  author  = {Agrachev, Andrei A. and Sarychev, Andrey V.},
  title   = {Controllability of 2D {Euler} and {Navier--Stokes} equations by degenerate forcing},
  journal = {Communications in Mathematical Physics},
  year    = {2006},
  volume  = {265},
  number  = {3},
  pages   = {673--697},
  doi     = {10.1007/s00220-006-0002-8}
}

@article{Shirikyan2007,
  author  = {Shirikyan, Armen},
  title   = {Exact controllability in projections for three-dimensional {Navier--Stokes} equations},
  journal = {Annales de l'Institut Henri Poincar{\'e} C, Analyse Non Lin{\'e}aire},
  year    = {2007},
  volume  = {24},
  number  = {4},
  pages   = {521--537},
  doi     = {10.1016/j.anihpc.2006.04.002}
}

@article{10540044,
  author  = {Adu, Daniel Owusu and Gharesifard, Bahman},
  title   = {Approximate controllability of continuity equation of {Transformers}},
  journal = {IEEE Control Systems Letters},
  year    = {2024},
  volume  = {8},
  pages   = {964--969},
  doi     = {10.1109/LCSYS.2024.3406255}
}

@article{Agrachev2025,
  author  = {Agrachev, Andrei and Letrouit, Cyril},
  title   = {Generic controllability of equivariant systems and applications to particle systems and neural networks},
  journal = {Annales de l'Institut Henri Poincar{\'e} C, Analyse non lin{\'e}aire},
  year    = {2026},
  volume  = {43},
  number  = {3},
  pages   = {639--668},
  doi     = {10.4171/AIHPC/151}
}

@book{agrachevControlTheoryGeometric2004,
  author    = {Agrachev, Andrei A. and Sachkov, Yuri L.},
  title     = {Control Theory from the Geometric Viewpoint},
  series    = {Encyclopaedia of Mathematical Sciences},
  volume    = {87},
  publisher = {Springer},
  address   = {Berlin, Heidelberg},
  year      = {2004},
  doi       = {10.1007/978-3-662-06404-7}
}

@article{agrachevControlManifoldsMappings2022,
  author  = {Agrachev, Andrei and Sarychev, Andrey},
  title   = {Control on the manifolds of mappings with a view to the deep learning},
  journal = {Journal of Dynamical and Control Systems},
  year    = {2022},
  volume  = {28},
  number  = {4},
  pages   = {989--1008},
  doi     = {10.1007/s10883-021-09561-2}
}

@article{agrachevControlSpacesEnsembles2020,
  author  = {Agrachev, Andrei and Sarychev, Andrey},
  title   = {Control in the spaces of ensembles of points},
  journal = {SIAM Journal on Control and Optimization},
  year    = {2020},
  volume  = {58},
  number  = {3},
  pages   = {1579--1596},
  doi     = {10.1137/19M1273049}
}

@article{AlcaldeFantuzziZuazua2026,
  author  = {Alcalde, Albert and Fantuzzi, Giovanni and Zuazua, Enrique},
  title   = {Exact sequence interpolation with {Transformers}},
  journal = {Mathematical Foundations of Machine Learning},
  year    = {2026},
  volume  = {2},
  number  = {1},
  pages   = {2},
  doi     = {10.1007/s44439-026-00005-y}
}

@book{AGS,
  author    = {Ambrosio, Luigi and Gigli, Nicola and Savar{\'e}, Giuseppe},
  title     = {Gradient Flows: In Metric Spaces and in the Space of Probability Measures},
  series    = {Lectures in Mathematics ETH Z{\"u}rich},
  publisher = {Birkh{\"a}user},
  address   = {Basel},
  year      = {2005},
  doi       = {10.1007/b137080}
}

@article{AVERBOUKH2025113205,
  author  = {Averboukh, Yurii and Khlopin, Dmitry},
  title   = {Pontryagin maximum principle for the deterministic mean field type optimal control problem via the {Lagrangian} approach},
  journal = {Journal of Differential Equations},
  year    = {2025},
  volume  = {430},
  pages   = {113205},
  doi     = {10.1016/j.jde.2025.02.076}
}

@article{BonnetFrankowska2021a,
  author  = {Bonnet, Beno{\^\i}t and Frankowska, H{\'e}l{\`e}ne},
  title   = {Necessary optimality conditions for optimal control problems in {Wasserstein} spaces},
  journal = {Applied Mathematics \& Optimization},
  year    = {2021},
  volume  = {84},
  number  = {S2},
  pages   = {1281--1330},
  doi     = {10.1007/s00245-021-09772-w}
}

@misc{bruno2024emergencemetastableclusteringmeanfield,
  author        = {Bruno, Giuseppe and Pasqualotto, Federico and Agazzi, Andrea},
  title         = {Emergence of meta-stable clustering in mean-field {Transformer} models},
  year          = {2024},
  eprint        = {2410.23228},
  archivePrefix = {arXiv},
  primaryClass  = {cs.LG},
  doi           = {10.48550/arXiv.2410.23228}
}

@book{CardMaster2019,
  author    = {Cardaliaguet, Pierre and Delarue, Fran{\c{c}}ois and Lasry, Jean-Michel and Lions, Pierre-Louis},
  title     = {The Master Equation and the Convergence Problem in Mean Field Games},
  series    = {Annals of Mathematics Studies},
  volume    = {201},
  publisher = {Princeton University Press},
  address   = {Princeton, NJ},
  year      = {2019},
  doi       = {10.23943/princeton/9780691190716.001.0001}
}

@book{CarmonaDelarue2018,
  author    = {Carmona, Ren{\'e} and Delarue, Fran{\c{c}}ois},
  title     = {Probabilistic Theory of Mean Field Games with Applications I: Mean Field FBSDEs, Control, and Games},
  series    = {Probability Theory and Stochastic Modelling},
  volume    = {83},
  publisher = {Springer},
  address   = {Cham},
  year      = {2018},
  doi       = {10.1007/978-3-319-58920-6}
}

@misc{CastinAblinCarrilloPeyre2025,
  author        = {Castin, Val{\'e}rie and Ablin, Pierre and Carrillo, Jos{\'e} Antonio and Peyr{\'e}, Gabriel},
  title         = {A unified perspective on the dynamics of deep {Transformers}},
  year          = {2025},
  eprint        = {2501.18322},
  archivePrefix = {arXiv},
  primaryClass  = {cs.LG},
  doi           = {10.48550/arXiv.2501.18322}
}

@article{cavagnariLagrangianEulerianKantorovich2022,
  author  = {Cavagnari, Giulia and Lisini, Stefano and Orrieri, Carlo and Savar{\'e}, Giuseppe},
  title   = {{Lagrangian}, {Eulerian} and {Kantorovich} formulations of multi-agent optimal control problems: Equivalence and {$\Gamma$}-convergence},
  journal = {Journal of Differential Equations},
  year    = {2022},
  volume  = {322},
  pages   = {268--364},
  doi     = {10.1016/j.jde.2022.03.019}
}

@article{CavagnariMarigonda2021Attainability,
  author  = {Cavagnari, Giulia and Marigonda, Antonio},
  title   = {Attainability property for a probabilistic target in {Wasserstein} spaces},
  journal = {Discrete and Continuous Dynamical Systems},
  year    = {2021},
  volume  = {41},
  number  = {2},
  pages   = {777--812},
  doi     = {10.3934/dcds.2020300}
}

@article{ChertovskihPogodaevStaritsynAguiar2023,
  author  = {Chertovskih, Roman and Pogodaev, Nikolay and Staritsyn, Maxim and Aguiar, A. Pedro},
  title   = {Optimal control of distributed ensembles with application to {Bloch} equations},
  journal = {IEEE Control Systems Letters},
  year    = {2023},
  volume  = {7},
  pages   = {2059--2064},
  doi     = {10.1109/LCSYS.2023.3284205}
}

@article{Chow1939,
  author  = {Chow, Wei-Liang},
  title   = {{\"U}ber Systeme von linearen partiellen Differentialgleichungen erster Ordnung},
  journal = {Mathematische Annalen},
  year    = {1939},
  volume  = {117},
  pages   = {98--105},
  doi     = {10.1007/BF01450011}
}

@article{Cuchiero,
  author  = {Cuchiero, Christa and Larsson, Martin and Teichmann, Josef},
  title   = {Deep neural networks, generic universal interpolation, and controlled {ODEs}},
  journal = {SIAM Journal on Mathematics of Data Science},
  year    = {2020},
  volume  = {2},
  number  = {3},
  pages   = {901--919},
  doi     = {10.1137/19M1284117}
}

@book{Deimling1977,
  author    = {Deimling, Klaus},
  title     = {Ordinary Differential Equations in Banach Spaces},
  series    = {Lecture Notes in Mathematics},
  volume    = {596},
  publisher = {Springer},
  address   = {Berlin, Heidelberg},
  year      = {1977},
  doi       = {10.1007/BFb0091636}
}

@book{dudley:89:real,
  author    = {Dudley, Richard M.},
  title     = {Real Analysis and Probability},

  publisher = {Wadsworth \& Brooks/Cole Advanced Books \& Software},
  address = {CA},
  year      = {1989}
}

@article{DudnikovSamborski1980,
  author  = {Dudnikov, P. I. and Samborski, S. N.},
  title   = {Criterion of controllability for systems in {Banach} space (generalization of {Chow}'s theorem)},
  journal = {Ukrain. Mat. Zh.},
  year    = {1980},
  volume  = {32},
  number  = {5},
  pages   = {649--653},
  note    = {In Russian}
}

@article{duprezApproximateExactControllability2019,
  author  = {Duprez, Michel and Morancey, Morgan and Rossi, Francesco},
  title   = {Approximate and exact controllability of the continuity equation with a localized vector field},
  journal = {SIAM Journal on Control and Optimization},
  year    = {2019},
  volume  = {57},
  number  = {2},
  pages   = {1284--1311},
  doi     = {10.1137/17M1152917}
}

@article{duprezMinimalTime2020,
  author  = {Duprez, Michel and Morancey, Morgan and Rossi, Francesco},
  title   = {Minimal time for the continuity equation controlled by a localized perturbation of the velocity vector field},
  journal = {Journal of Differential Equations},
  year    = {2020},
  volume  = {269},
  number  = {1},
  pages   = {82--124},
  doi     = {10.1016/j.jde.2019.11.098}
}

@inproceedings{FuchsWorrallFischerWelling2020,
  author    = {Fuchs, Fabian B. and Worrall, Daniel E. and Fischer, Volker and Welling, Max},
  title     = {{SE(3)}-Transformers: {3D} roto-translation equivariant attention networks},
  booktitle = {Advances in Neural Information Processing Systems},
  volume    = {33},
  year      = {2020}
}

@article{FuruyaDeHoopLassas2025,
  author  = {Furuya, Takashi and {de Hoop}, Maarten V. and Lassas, Matti},
  title   = {{Transformers} through the lens of support-preserving maps between measures},
  journal = {arXiv preprint arXiv:2509.25611},
  year    = {2025},
  doi     = {10.48550/arXiv.2509.25611}
}

@article{Geshkovski2024,
  author  = {Geshkovski, Borjan and Letrouit, Cyril and Polyanskiy, Yury and Rigollet, Philippe},
  title   = {A mathematical perspective on {Transformers}},
  journal = {Bulletin of the American Mathematical Society},
  year    = {2025},
  volume  = {62},
  number  = {3},
  pages   = {427--479},
  doi     = {10.1090/bull/1863}
}

@article{GeshkovskiRigolletRuizBalet2026,
  author  = {Geshkovski, Borjan and Rigollet, Philippe and Ruiz-Balet, Dom{\`e}nec},
  title   = {Measure-to-measure interpolation using {Transformers}},
  journal = {Foundations of Computational Mathematics},
  year    = {2026},
  doi     = {10.1007/s10208-026-09754-w}
}

@article{Hermes1979,
  author  = {Hermes, Henry},
  title   = {Local controllability of observables in finite and infinite dimensional nonlinear control systems},
  journal = {Applied Mathematics and Optimization},
  year    = {1979},
  volume  = {5},
  pages   = {117--125},
  doi     = {10.1007/BF01442549}
}

@article{JimenezMarigondaQuincampoix2020,
  author  = {Jimenez, Chlo{\'e} and Marigonda, Antonio and Quincampoix, Marc},
  title   = {Optimal control of multiagent systems in the {Wasserstein} space},
  journal = {Calculus of Variations and Partial Differential Equations},
  year    = {2020},
  volume  = {59},
  number  = {2},
  pages   = {58},
  doi     = {10.1007/s00526-020-1718-6}
}

@inproceedings{Kaba2023,
  author    = {Kaba, S{\'e}kou-Oumar and Mondal, Arnab Kumar and Zhang, Yan and Bengio, Yoshua and Ravanbakhsh, Siamak},
  title     = {Equivariance with learned canonicalization functions},
  booktitle = {Proceedings of the 40th International Conference on Machine Learning},
  series    = {Proceedings of Machine Learning Research},
  volume    = {202},
  pages     = {15546--15566},
  year      = {2023}
}

@article{KhajehSalehaniMarkina2014,
  author  = {Khajeh Salehani, Mahdi and Markina, Irina},
  title   = {Controllability on infinite-dimensional manifolds: A {Chow--Rashevsky} theorem},
  journal = {Acta Applicandae Mathematicae},
  year    = {2014},
  volume  = {134},
  pages   = {229--246},
  doi     = {10.1007/s10440-014-9880-5}
}

@article{KipkaLedyaev2015,
  author  = {Kipka, Robert J. and Ledyaev, Yuri S.},
  title   = {Extension of chronological calculus for dynamical systems on manifolds},
  journal = {Journal of Differential Equations},
  year    = {2015},
  volume  = {258},
  number  = {5},
  pages   = {1765--1790},
  doi     = {10.1016/j.jde.2014.11.014}
}

@inproceedings{KohlerKleinNoe2020,
  author    = {K{\"o}hler, Jonas and Klein, Leon and No{\'e}, Frank},
  title     = {Equivariant flows: Exact likelihood generative learning for symmetric densities},
  booktitle = {Proceedings of the 37th International Conference on Machine Learning},
  series    = {Proceedings of Machine Learning Research},
  volume    = {119},
  pages     = {5361--5370},
  year      = {2020}
}

@book{kolokoltsov2019,
  author    = {Kolokoltsov, Vassili N.},
  title     = {Differential Equations on Measures and Functional Spaces},
  series    = {Birkh{\"a}user Advanced Texts: Basler Lehrb{\"u}cher},
  publisher = {Birkh{\"a}user},
  address   = {Cham},
  year      = {2019},
  doi       = {10.1007/978-3-030-03377-4}
}

@book{Lang1999,
  author    = {Lang, Serge},
  title     = {Fundamentals of Differential Geometry},
  series    = {Graduate Texts in Mathematics},
  volume    = {191},
  publisher = {Springer},
  address   = {New York},
  year      = {1999},
  doi       = {10.1007/978-1-4612-0541-8}
}

@article{Ledyaev2004,
  author  = {Ledyaev, Yuri S.},
  title   = {On an infinite-dimensional variant of the {Chow--Rashevskii} theorem},
  journal = {Doklady Mathematics},
  year    = {2004},
  volume  = {70},
  number  = {2},
  pages   = {782--784}
}

@article{lottTangentConesWasserstein2016,
  author  = {Lott, John},
  title   = {On tangent cones in {Wasserstein} space},
  journal = {Proceedings of the American Mathematical Society},
  year    = {2017},
  volume  = {145},
  number  = {7},
  pages   = {3127--3136},
  doi     = {10.1090/proc/13415}
}

@article{Mauhart1992CommutatorsOF,
  author  = {Mauhart, Markus and Michor, Peter W.},
  title   = {Commutators of flows and fields},
  journal = {Archivum Mathematicum},
  year    = {1992},
  volume  = {28},
  number  = {3--4},
  pages   = {229--236}
}

@article{ottoGeometryDissipativeEvolution2001,
  author  = {Otto, Felix},
  title   = {The geometry of dissipative evolution equations: The porous medium equation},
  journal = {Communications in Partial Differential Equations},
  year    = {2001},
  volume  = {26},
  number  = {1--2},
  pages   = {101--174},
  doi     = {10.1081/PDE-100002243}
}

@article{PogodaevRossi2024TrajectoryStabilization,
  author  = {Pogodaev, Nikolay and Rossi, Francesco},
  title   = {Trajectory stabilization of nonlocal continuity equations by localized controls},
  journal = {SIAM Journal on Control and Optimization},
  year    = {2024},
  volume  = {62},
  number  = {6},
  pages   = {3315--3340},
  doi     = {10.1137/24M1644274}
}

@article{PogodaevRossi2026GeodesicStabilization,
  author  = {Pogodaev, Nikolay and Rossi, Francesco},
  title   = {Stabilization of nonlocal continuity equations via localized controls: Case of geodesic perturbations},
  journal = {Systems \& Control Letters},
  year    = {2026},
  volume  = {211},
  pages   = {106396},
  doi     = {10.1016/j.sysconle.2026.106396}
}

@inproceedings{Puny2022,
  author    = {Puny, Omri and Atzmon, Matan and Smith, Edward J. and Misra, Ishan and Grover, Aditya and Ben-Hamu, Heli and Lipman, Yaron},
  title     = {Frame averaging for invariant and equivariant network design},
  booktitle = {10th International Conference on Learning Representations (ICLR)},
  year      = {2022}
}

@article{Rashevskii1938,
  author  = {Rashevsky, P. K.},
  title   = {Any two points of a totally nonholonomic space may be connected by an admissible line},
  journal = {Uch. Zap. Mosk. Gos. Ped. Inst. im. K. Libknekhta, Ser. Fiz.-Mat. Nauk},
  year    = {1938},
  volume  = {2},
  pages   = {83--94},
  note    = {In Russian}
}

@article{StaritsynPogodaevChertovskihPereira2021,
  author  = {Staritsyn, Maxim and Pogodaev, Nikolay and Chertovskih, Roman and Pereira, Fernando Lobo},
  title   = {Feedback maximum principle for ensemble control of local continuity equations: An application to supervised machine learning},
  journal = {IEEE Control Systems Letters},
  year    = {2022},
  volume  = {6},
  pages   = {1046--1051},
  doi     = {10.1109/LCSYS.2021.3089139}
}

@article{STEPANOV20171044,
  author  = {Stepanov, Eugene and Trevisan, Dario},
  title   = {Three superposition principles: Currents, continuity equations and curves of measures},
  journal = {Journal of Functional Analysis},
  year    = {2017},
  volume  = {272},
  number  = {3},
  pages   = {1044--1103},
  doi     = {10.1016/j.jfa.2016.10.025}
}

@inproceedings{Vaswani2017,
  author    = {Vaswani, Ashish and Shazeer, Noam and Parmar, Niki and Uszkoreit, Jakob and Jones, Llion and Gomez, Aidan N. and Kaiser, {\L}ukasz and Polosukhin, Illia},
  title     = {Attention is all you need},
  booktitle = {Advances in Neural Information Processing Systems},
  volume    = {30},
  pages     = {5998--6008},
  year      = {2017}
}

\end{document}